\documentclass[12pt,twoside,english]{article}
\usepackage{a4wide}
\pdfoutput=1

\makeatletter
\usepackage[utf8]{inputenc}            
\usepackage[T1]{fontenc} 			
\usepackage{graphicx}
\usepackage{fancyhdr}
\usepackage{amsmath, amssymb, amsfonts, bm}
\usepackage{mathtools}
\usepackage{caption}
\usepackage[labelformat=simple]{subfig}
\usepackage{setspace}
\usepackage{verbatim} % enables multi-line comments.
\usepackage[english]{babel}
\usepackage[novbox]{pdfsync} %may cause large empty boxes
\usepackage{caption}
\usepackage{booktabs}
\usepackage{hyperref}
\usepackage{float}

\usepackage{siunitx}
\usepackage{xfrac}
\usepackage{pdfpages}

\usepackage{color}

\includepdfset{turn=true,pagecommand={}}

\hypersetup{
	colorlinks = false,
	linktocpage = true,
	allbordercolors = {0.8 0.8 0.8}
}

\addto\captionsenglish{}
\addto\captionsenglish{}
\addto\extrasenglish{}
\addto\extrasenglish{}
\addto\extrasenglish{}
\addto\extrasenglish{}
\addto\extrasenglish{}
\addto\extrasenglish{}
\addto\extrasenglish{}

\usepackage{etex}
\usepackage{tikz}
\usepackage{pgfplots}
\usepackage{tikz-3dplot}
\usepackage{fp}

\pgfplotsset{compat=newest}

\ifdefined\bUseTikzExternalize
	\usetikzlibrary{external}
	\tikzsetexternalprefix{tmp/}
\else
	\usetikzlibrary{external}
	\tikzexternaldisable
\fi

\usetikzlibrary{calc,intersections}
\usetikzlibrary{shapes.geometric,shapes.arrows,decorations.pathmorphing}
\usetikzlibrary{matrix,chains,scopes,positioning,arrows,fit}
\usetikzlibrary{spy}
\usetikzlibrary{fadings}
\usetikzlibrary{shadings}
\usetikzlibrary{backgrounds}

\usepgfplotslibrary{groupplots}
\usepgfplotslibrary{patchplots}

\def\pgfplotfontsizetitle{\small}
\def\pgfplotfontsize{\small}

\pgfplotsset{
  mystyle/.style ={%
    grid = major,
    every tick label/.append style={font=\pgfplotfontsize},
    every axis label/.append style={font=\pgfplotfontsize},
    legend style={font=\scriptsize},
    label style={font=\pgfplotfontsize},
    title style={font=\pgfplotfontsizetitle},
    /pgf/number format/set thousands separator = {}, % 1,000 -> 1000
  }
}

\tikzset{
	>=stealth',
	axis/.style={<->},
	important line/.style={thick},
	connection/.style={thick, dotted},
	dot/.style = {
		draw,
		fill = white,
		circle,
		solid,
		thin,
		inner sep = 0pt,
		minimum size = 3pt
	},
	defP/.style = {
		inner sep = 0pt
	}	
}
\pgfarrowsdeclaredouble{doubled}{doubled}{stealth'}{stealth'}

\usepackage[outline]{contour}
\usepackage{standalone} % for standalone compilation of tikz graphics
\newcommand{\vek}[1]{\mathchoice{\displaystyle\boldsymbol#1}
	{\textstyle\boldsymbol#1}{\scriptstyle\boldsymbol#1}
	{\scriptscriptstyle\boldsymbol#1}}
\newcommand{\mat}[1]{\mathchoice{\displaystyle\mathbf#1}
	{\textstyle\mathbf#1}{\scriptstyle\mathbf#1}
	{\scriptscriptstyle\mathbf#1}}
\newcommand{\ten}[1]{ \ensuremath{\bm #1} }
\renewcommand{\d}{ \ensuremath{\mathrm{d} }}
\newcommand{\p}{ \ensuremath{\partial} }
\renewcommand{\t}[1]{ \ensuremath{\text{#1}} }
\newcommand{\T}{{\ensuremath{\mathrm{T}}}}

\newcommand{\divG}[1]{\ensuremath{\text{div}_{\Gamma} #1}}
\newcommand{\gradR}[1]{\ensuremath{\nabla_{\vek{r}} #1}}
\newcommand{\gradG}[1]{\ensuremath{\nabla_\Gamma #1}}
\newcommand{\gradGD}[1]{\ensuremath{\nabla_\Gamma^\text{dir} #1}}
\newcommand{\gradGC}[1]{\ensuremath{\nabla_\Gamma^\text{cov} #1}}
\newcommand{\nG}{\ensuremath{\vek{n}_\Gamma}}
\newcommand{\nCo}{\ensuremath{\vek{n}_{\p\Gamma}}}
\newcommand{\tB}{\ensuremath{\vek{t}_{\p\Gamma}}}

\newcommand{\vu}{\ensuremath{\vek{v}_u}}
\newcommand{\vw}{\ensuremath{\vek{v}_w}}

\newcommand{\Nu}{\ensuremath{\vek{N}_{\vek{u}}}}
\newcommand{\Nw}{\ensuremath{\vek{N}_{\vek{w}}}}

\newcommand{\lyxaddress}[1]{
	\par {\raggedright #1
		\vspace{1.4em}
		\noindent\par}
}

\numberwithin{equation}{section}

\makeatother

\graphicspath{{Pics/}}
\def\pathToBibFile{References}

\newcommand{\revStart}{\color{black}}
\newcommand{\revEnd}{\color{black}}

\begin{document}

\title{Reissner--Mindlin shell theory based on tangential differential calculus}

\author{D. Sch{\"o}llhammer, T.P. Fries}
\maketitle

\lyxaddress{\begin{center}
Institute of Structural Analysis\\
Graz University of Technology\\
Lessingstr. 25/II, 8010 Graz, Austria\\
\texttt{www.ifb.tugraz.at}\\
\texttt{schoellhammer@tugraz.at}
\end{center}}

\begin{abstract}

The linear Reissner-Mindlin shell theory is reformulated in the frame of the tangential differential calculus (TDC) using a global Cartesian coordinate system. The rotation of the normal vector is modelled with a difference vector approach. The resulting equations are applicable to both explicitly and implicitly defined shells, because the employed surface operators do not necessarily rely on a parametrization. Hence, shell analysis on surfaces implied by level-set functions is enabled, but also the classical case of parametrized surfaces is captured. As a consequence, the proposed TDC-based formulation is more general and may also be used in recent finite element approaches such as the TraceFEM and CutFEM where a parametrization of the middle surface is not required. Herein, the numerical results are obtained by isogeometric analysis using NURBS as trial and test functions for classical and new benchmark tests. In the residual errors, optimal higher-order convergence rates are confirmed when the involved physical fields are sufficiently smooth.


\end{abstract} 
\newpage\tableofcontents\newpage
\section{Introduction}
\label{sec:intro}

Shells are curved, thin-walled structures that occur in a wide range of applications in both nature and technology. Because they often feature a high bearing capacity, they are frequently used in engineering applications such as automotive, aerospace, biomedical- and civil-engineering, see e.g., \cite{Calladine_1983a,Farshad_1992a,Zingoni_2018a}. For the mechanical modelling, one may distinguish two classes of models based on the fact whether traverse shear deformations are considered or not. The Kirchhoff-Love shell theory \cite{Kirchhoff_1850a,Love_1888a} neglects such shear deformation resulting in a fourth-order partial differential equation (PDE) for the displacement of the middle surface of thin shells. For the numerical treatment, it is crucial to consider continuity requirements due to the variational index 2 \cite{Echter_2013a}. When taking shear deformations into account, the Reissner-Mindlin shell theory \cite{Reissner_1945a} is typically employed and models the deformation of thin \emph{and} moderately thick shells. The resulting shell equations are a system of second-order PDEs with the unknowns being the displacement of the middle surface and the rotation of the normal vector. An advantage of this model is that the corresponding variational index is $1$ and only $C^0$-continuity in a finite element analysis is required. On the other hand, the model often suffers from locking phenomena for increasingly thin shells which may require  further measures in the numerical treatment. An overview of classical shell theory is given, e.g., in \cite{Chapelle_2000a,Simo_1989a,Simo_1989b,Bischoff_2017a} or in the text books \cite{Calladine_1983a,Blaauwendraad_2014a,Basar_1985a,Wempner_2002a,Zingoni_2018a}.\par

The middle surface of a shell is a manifold of codimension 1 embedded in the physical space $\mathbb{R}^3$. One may distinguish two alternatives for the definition of the shells: (1) explicit or (2) implicit. The first approach is typically used in the classical shell theory, where the middle surface of the shell is given through a parametrization, i.e., a map from a two-dimensional parameter space to the real surface embedded in the physical space $\mathbb{R}^3$. Based on this parametrization, curvilinear coordinates may be introduced and geometric quantities (normal vectors, curvatures, etc.) and differential surface operators (gradients, divergence, etc.) are defined. These quantities are the key ingredients of the PDEs obtained in the classical shell theory. There are various implementations of the classical shell models using this approach, a general overview is given, e.g., in \cite{Yang_2000a}, Kirchhoff-Love shells may be found, e.g., in \cite{Cirak_2000a,Cirak_2001a,Kiendl_2009a,Nguyen-Thanh_2017a,Tepole_2015a}, Reissner-Mindlin shells, e.g., in \cite{Benson_2010a,Dornisch_2013a,Dornisch_2014a,Bathe_2000a,Ko_2017a,Kiendl_2017a} and hierarchical shells, e.g., in \cite{Echter_2013a,Oesterle_2016a,Bieber_2018a}.\par

The alternative to parametrizations is when the middle surface is defined \emph{implicitly}, e.g.~as the zero-isosurface of a scalar function in three dimensions following the level-set method \cite{Fries_2017a,Fries_2017b,Osher_2003a,Sethian_1999b}. The geometric quantities and differential surface operators are defined in the global Cartesian coordinates system of the surrounding physical space as done in the tangential differential calculus (TDC) \cite{Delfour_2011a,Gurtin_1975a,Hansbo_2015a,Schoellhammer_2018a}. When the physical modeling process is based on the TDC, it is applicable to surfaces which are parametrized (explicit definition) or \emph{not} (implicit definition). In this sense, the TDC-based approach is more general than approaches based on local coordinates, which are restricted to explicit surface descriptions. Models based on the TDC are found in various applications, see \cite{Demlow_2009a,Dziuk_1988a,Dziuk_2013a,Fries_2017b} for scalar problems such as heat flow and \cite{Fries_2018b,Jankuhn_2017a} for flow problems on manifolds. In the context of structure mechanics, this approach is used in \cite{Delfour_1994a,Delfour_1995a,Delfour_1996a,Opstal_2015a,Hansbo_2017a,Schoellhammer_2018a} for Kirchhoff-Love shells, in \cite{Hansbo_2014b} for curved beams, and in \cite{Hansbo_2014a,Hansbo_2015a,Gurtin_1975a} for membranes. In addition to a more general formulation of the PDEs in the frame of the TDC, there is unified and elegant implementation in the finite element context that may be recycled in other situations where PDEs on surfaces are considered. More precisely, typical surface operators in the TDC, in particular, the surface gradient may be applied to the finite element shape functions independently of the concrete application. This enables a shift of significant parts of the implementation needed for shells to the underlying finite element technology.\par

In \cite{Schoellhammer_2018a}, the authors proposed a reformulation of the Kirchhoff-Love shell based on the TDC. Herein, the aim is to recast the Reissner-Mindlin shell in the TDC-framework using a difference vector formulation. All relevant mechanical fields such as membrane forces, bending moments and shear forces are expressed in the global Cartesian coordinate system using the TDC and the computation of invariant quantities such as principal moments is shown. Furthermore, a parametrization-free strong form of the force and moment equilibriums is derived and used as a starting point for the derivation of the weak form. Similar to \cite{Schoellhammer_2018a}, the concept of residual errors is also used to confirm higher-order convergence rates in the corresponding error norm for suitable shell test cases. This is, otherwise, very difficult as exact solutions for shells are scarce and standard benchmarks hardly serve for higher-order convergence studies.\par

For the numerical solution of the shell boundary value problem (BVP), one may employ two fundamentally different approaches of solving PDEs on manifolds. The first approach is a classical finite element method, labelled Surface FEM herein \cite{Demlow_2009a,Dziuk_2013a,Fries_2018b,Fries_2017b}. In the classical Surface FEM the analytical geometry is decomposed into a set of elements. Each element of the discrete surface implies an element-wise parametrization, no matter if the geometry was originally defined explicitly or implicitly. In this case, the classical shell theory based on local coordinates \emph{and} the TDC-based formulation are suitable choices provided that the geometry is discretized by finite elements. The alternative numerical approach is to embed the implicitly defined shell surface in a three-dimensional background mesh. Boundaries of the shell can be defined by means of the boundary of the background mesh or with additional level-set functions as in \cite{Fries_2017a,Fries_2017b}. The trial and test functions are the three-dimensional shape functions of the background mesh restricted to the trace of the shell surface. These methods are labelled as CutFEM \cite{Burman_2015a,Burman_2018a,Cenanovic_2016a,Elfverson_2018a} or TraceFEM \cite{Grande_2016a,Olshanskii_2017a,Olshanskii_2017b,Reusken_2014a}. When applying these methods to shell mechanics, it is crucial to use a TDC-based formulation such as proposed in this work, because no parametrization is available nor needed.\par

The presented numerical results herein are obtained with the Surface FEM \cite{Demlow_2009a,Dziuk_2013a,Fries_2018b,Fries_2017b} using NURBS as trial and test functions as proposed by Hughes et al.~\cite{Cottrell_2009a,Hughes_2005a} in a general sense. In \cite{Benson_2010a,Dornisch_2013a,Dornisch_2014a,Kiendl_2017a}, isogeometric analysis (IGA) is also applied to the Reissner-Mindlin shell, however, based on the classical, parameter-based formulation of the governing equations. \revStart The same results are obtained here, however, based on the new TDC-based formulation resulting into a significant different implementation. In particular, the implementation is more intuitive and compact due to a sharper split of FE-technology and application, as already shown in \cite{Schoellhammer_2018a}.\revEnd\ Although IGA is used in the numerical results shown here, we emphasize that the presented TDC-based formulation, being the main aspect of this work, is more general as it may also be used when no parametrization is available such as in CutFEM or TraceFEM.\par

It is pointed out that continuity requirements would also allow a standard FEM implementation using $C^0$-continuous shape functions. Nevertheless, we prefer the use of NURBS here, for example, due to the continuous normal vector and improved convergence properties. As mentioned above, a standard difference vector formulation is employed here and different strategies of discretizing tangential vector fields are outlined. We use Lagrange multipliers in order to weakly enforce (i) the tangentiality constraint on the globally defined difference vector and (ii) the boundary conditions. The proposed approach is in the case of parametrized surfaces equivalent to the classical Reissner-Mindlin shell. Therefore, locking phenomena can be expected in the case of thin shells. As shown in the numerical results, order elevation can quite efficiently reduce the locking effects and is easily achieved with the isogeometric approach. Further treatment of locking is not considered herein and is beyond of the scope of this work.\par

The outline of the paper is as follows: In \autoref{sec:pre}, a brief introduction to the tangential differential calculus (TDC) is given. In \autoref{sec:method}, the classical linear Reissner-Mindlin shell equations are recast in terms of the TDC. Stress resultants such as membrane forces, bending moments and transverse shear forces are defined. Based on the stress resultants, the force and moment equilibriums are presented. In \autoref{sec:impl}, the discrete weak form and the resulting system of linear equations is shown. Furthermore, different discetization strategies of the difference vector are outlined. In \autoref{sec:numres}, numerical results are presented. The first example is the well-known Scordelis-Lo roof proposed in \cite{Belytschko_1985a}, the second example is the partly clamped hyperbolic paraboloid taken from \cite{Bathe_2000a}, and the last example is based on the clamped flower shaped shell from \cite{Schoellhammer_2018a}. The error is measured in the strong form of the equilibrium in order to verify the proposed approach and higher-order convergence rates are achieved.

\section{Preliminaries}
\label{sec:pre}

\begin{figure}[h]
	\centering
	\subfloat[explicit]{
	% file=#1, size=#2, figno=#3
			\ifdefined\bUseTikzExternalize
				\def\tkzscale{0.2}
				\centering
				\tikzsetnextfilename{Fig1a}  
				\input{tikz/Fig1a}
			\else
				\includegraphics[width=0.37\textwidth]{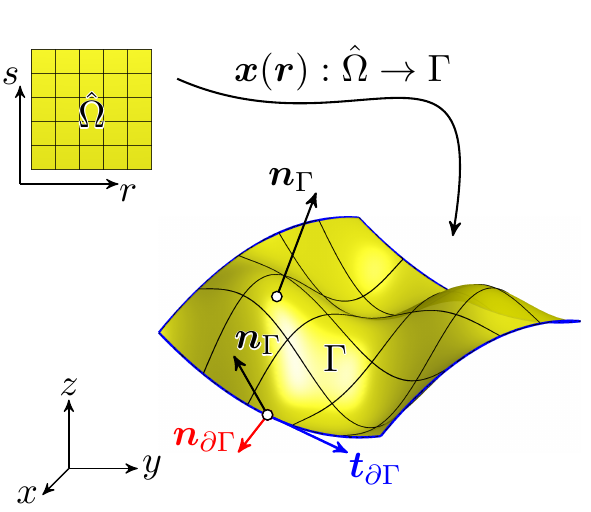}
			\fi			
 \label{fig:surfexpl}}
	\hfil
	\subfloat[implicit]{
	% file=#1, size=#2, figno=#3
			\ifdefined\bUseTikzExternalize
				\def\tkzscale{0.2}
				\centering
				\tikzsetnextfilename{Fig1b}  
				\input{tikz/Fig1b}
			\else
				\includegraphics[width=0.45\textwidth]{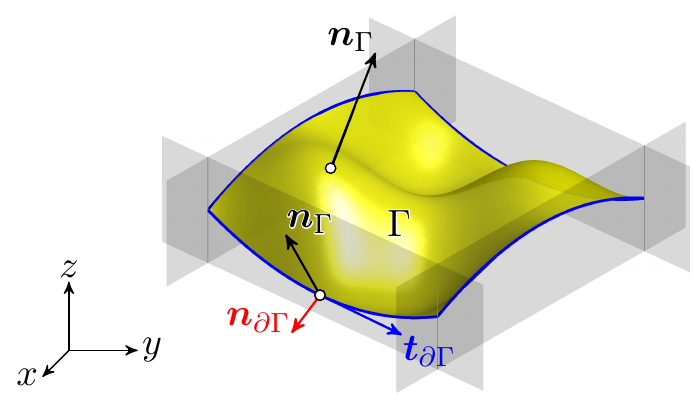}
			\fi			
\label{fig:surfimpl}}
	\caption{Examples of bounded surfaces $\Gamma$ embedded in the physical space $\mathbb{R}^3$: (a) explicitly defined surface with a parametrization $\vek{x}(\vek{r})$, (b) implicitly defined surface with a master level-set function $\phi(\vek{x}) = 0$ (yellow) and slave level-set functions $\psi_i$ for the boundary definition (grey).}
	\label{fig:surf}
\end{figure}

A shell is a thin-walled structure, which can be modelled as a surface $\Gamma$ embedded in the physical space $\mathbb{R}^3$. Let the middle surface of the shell be possibly curved, sufficiently smooth, orientable, connected and bounded by $\p\Gamma$. The surface may also be called a manifold. Surfaces may be defined (explicitly) through a bijective mapping $\vek{x}(\vek{r}) : \hat{\Omega} \to \Gamma$ from the parameter space $\hat{\Omega} \subset \mathbb{R}^2$ to the real domain $\Gamma \subset \mathbb{R}^3$, see \autoref{fig:surfexpl}. \revStart In the case of a surface mesh, this is rather an atlas of local, element-wise mappings. In these cases, \revEnd the surface is given by a parametrization or \emph{parametrized}. A surface may also be defined \emph{implicitly} by a level-set function $\phi(\vek{x}) : \mathbb{R}^3 \to \mathbb{R}$ following the level-set method. Then, $\phi$ is a scalar-valued function and the middle surface is implied by its zero-isosurface $\phi(\vek{x}) = 0$, which might be bounded by additional level-set functions as described in \cite{Fries_2017b}, see \autoref{fig:surfimpl}.  For the formal proof of equivalence of both cases we refer to, e.g., \cite{Dziuk_2013a}.\par
In both cases, there exists a unit normal vector $\vek{n} = [n_x,\,n_y,\,n_z]^\T \in \mathbb{R}^3$. The \revStart representation \revEnd of the normal vector depends on whether the surface is based on a parametrization or not. In the parametrized case, the normal vector is determined by the normalized cross-product of the tangent vectors living in the shell surface and given by the columns of the Jacobi-matrix $\mat{J}(\vek{r}) = \sfrac{\p\vek{x}}{\p\vek{r}}$. In the implicit case, the normal vector is \revStart given \revEnd by the normalized gradient of the level-set function $\nG = \sfrac{\nabla\phi}{\Vert\nabla\phi\Vert}$. Along the boundary $\p\Gamma$, there is an associated tangential vector $\vek{t}_{\p\Gamma} \in \mathbb{R}^3$ pointing along $\p\Gamma$ and a co-normal vector $\vek{n}_{\p\Gamma} =\nG \times \vek{t}_{\p\Gamma} \in \mathbb{R}^3$ pointing ``outwards'' and being perpendicular to the boundary $\p\Gamma$ yet in the tangent plane of the surface $\Gamma$.\par

\subsection{Surface gradients and divergence\label{sub:SurfaceOperators}}

On the manifold $\Gamma$, the orthogonal projection operator $\mat P\left(\vek x\right)\in\mathbb{R}^{3\times3}$
is defined based on the normal vector as 
\begin{align}
\mat P\left(\vek x\right)=\mathbb{I} - \nG\left(\vek x\right)\otimes\nG\left(\vek x\right)
\end{align}
and projects arbitrary vectors $\vek{v}(\Gamma) \in \mathbb{R}^3$ into the tangent space $T_P\Gamma$. There holds $\mat P\cdot\vek n_{\Gamma}=\vek0$, $\mat P=\mat P^{\T}$, and $\mat P\cdot\mat P=\mat P$. Instead, the projection operator
\begin{align}
\mat{Q}\left(\vek x\right) = \mathbb{I} - \mat P(\vek x)
\end{align}
projects an arbitrary vector in normal direction of $\Gamma$ and there holds $\mat{Q} \cdot \nG = \nG$, $\mat{Q} = \mat{Q}^\T$, and $\mat{Q}\cdot\mat{Q} = \mat{Q}$.\par

The tangential gradient operator $\nabla_{\Gamma}$ of a differentiable
\emph{scalar} function $u:\Gamma\to\mathbb{R}$ on the manifold is
given by
\begin{align}
\gradG{u(\vek x)} & = \mat P\left(\vek x\right)\cdot\nabla\tilde{u}\left(\vek x\right)\ , \quad \vek{x} \in \Gamma\label{eq:TangGradImplicit}
\end{align}
where $\nabla$ is the standard gradient operator, and $\tilde{u}$
is a smooth extension of $u$ in a neighborhood $\mathcal{U}$ of
the manifold $\Gamma$. For parametrized surfaces defined by the map
$\vek x\left(\vek r\right)$, and
a given scalar function $u\left(\vek r\right):\mathbb{R}^{2}\to\mathbb{R}$,
the tangential gradient may be determined without explicitly computing
an extension $\tilde{u}$ using
\begin{align}
\gradG{u(\vek x(\vek r))}=\mat J\left(\vek r\right)\cdot\mat G^{-1}\left(\vek r\right)\cdot\gradR{u(\vek r)}\ ,\label{eq:TangGradExplicit}
\end{align}
with $\mat G=\mat J^\T\cdot\mat J$ being
the metric tensor (first fundamental form). This is relevant in the
context of the Surface FEM, where $u\left(\vek r\right)$ may be the
shape functions in the reference element and tangential gradients
are to be determined in the physical surface elements. It is noteworthy
that $\gradG{u}$ is in the tangent space of $\Gamma$, i.e., $\gradG{u} \in 
T_P\Gamma$, and,
thus, $\mat P\cdot\nabla_{\Gamma}u=\gradG{u}$, $ \mat{Q} \cdot \gradG{u} = \vek{0}$ and $\gradG{u}\cdot\vek n_{\Gamma}=0$.
It is straightforward to determine second order derivatives of scalar
functions, see, e.g., \cite{Delfour_2011a,Schoellhammer_2018a}.\par

When surface gradients of \emph{vector} functions $\vek u\left(\vek x\right):\Gamma\to\mathbb{R}^{3}$
are considered, it is important to distinguish directional and covariant
gradients:
\begin{align*}	\gradGD{\vek{u}\left(\vek x\right)} &= \gradGD{\begin{bmatrix}
		u\left(\vek x\right)\\
		v\left(\vek x\right)\\
		w\left(\vek x\right)
	\end{bmatrix}}=\begin{bmatrix}
	\left(\gradG{u}\right)^\T\\
	\left(\gradG{v}\right)^\T\\
	\left(\gradG{w}\right)^\T
\end{bmatrix}=\nabla\tilde{\vek u}\cdot\mat P,\\
\nabla_{\Gamma}^{\mathrm{cov}}\vek u\left(\vek x\right) &= \mat P\cdot\gradGD{\vek u\left(\vek x\right)}=\mat P\cdot\nabla\tilde{\vek u}\cdot\mat P.
\end{align*}
Concerning the surface divergence of vector functions $\vek u\left(\vek x\right):\Gamma\to\mathbb{R}^{3}$
and tensor functions $\mat A\left(\vek x\right):\Gamma\to\mathbb{R}^{3\times3}$,
there holds 
\begin{align*}
	\divG{\vek u\left(\vek x\right)} & =  \divG{\left(u,v,w\right)^\T} = \mathrm{tr}\left(\gradGD{
	\vek u}\right)=\mathrm{tr}\left(\gradGC{\vek u}\right) \eqqcolon\nabla_{\Gamma}\cdot\vek u,\\
	\divG{\mat A\left(\vek x\right)} & =  \begin{bmatrix}
		\divG{\left(A_{11},A_{12},A_{13}\right)}\\
		\divG{\left(A_{21},A_{22},A_{23}\right)}\\
		\divG{\left(A_{31},A_{32},A_{33}\right)}
	\end{bmatrix}\eqqcolon\nabla_{\Gamma}\cdot\mat A.
\end{align*}
To derive the weak form of the governing equations, the following
divergence theorem on manifolds is needed \cite{Delfour_1996a,Delfour_2011a,Schoellhammer_2018a},
\begin{align}
\int_{\Gamma}\vek u\cdot\divG{\mat A}\, \d A=-\int_{\Gamma}\gradGD{\vek u}:\mat A\, \d A + \int_{\Gamma}\varkappa\cdot\vek u\cdot\mat A\cdot \nG\,\mathrm{d}A+\int_{\p\Gamma}\vek u\cdot\mat A\cdot \nCo\,\mathrm{d}s,\label{eq:DivTheorem}
\end{align}
where $\gradGD{\vek u}:\mat A=\mathrm{tr}\left(\gradGD{\vek u}\cdot\mat A^{\T}\right)$.
For \emph{tangential} (in-plane) tensor functions with $\mat A=\mat P\cdot\mat A\cdot\mat P$,
the term involving the mean curvature $\varkappa = \t{tr}(\gradGD{\nG})$ vanishes and one finds
$\gradGD{\vek u}:\mat A=\gradGC{\vek u}:\mat A$.

\section{The shell equations}
\label{sec:method}

In this section, we derive the linear Reissner-Mindlin shell theory or first-order shear deformation theory in the frame of tangential operators based on a global Cartesian coordinate system. As mentioned before, this has the advantage over classical shell theory that the resulting model is valid no matter whether a parametrization is available or not.\par
We restrict ourselves to infinitesimal deformations and rotations, which means that the reference and spatial configuration are indistinguishable. For simplicity, a linear elastic material governed by Hooke's law is assumed. In contrast to the Kirchhoff-Love shells, the additional constraint on the shell director is omitted, thus allowing transverse shear strains, leading to the well-known 5-parameter shell models, e.g., \cite{Bischoff_2017a}. Furthermore, we assume a constant shifter in the material law, which enables an analytical pre-integration in thickness direction.\par

The shell continuum $\Omega$ of thickness $t$ can be defined implicitly \revStart by a special level-set function, called signed-distance function $\phi_{\t{SDF}}$. This scalar-valued function gives the shortest Euclidean distance from $\vek{x}$ to the middle surface of the shell and is positive in direction of the normal vector $\nG$ and negative otherwise. Then,
\begin{align}
\Omega = \left\lbrace \vek{x} \in \mathbb{R}^3 : \vert \phi_{\t{SDF}}(\vek{x}) \vert \le \frac{t}{2} \right\rbrace\ .
\end{align}
\revEnd Alternatively, when the middle surface $\Gamma$ is parametrized with a map $\vek{x}_\Gamma(\vek{r})$, the domain of the shell is defined by
\begin{align}\label{eq:xExpl}
\vek{x} = \vek{x}_\Gamma + \zeta \nG(\vek{x}_\Gamma)\ ,
\end{align}
where $\zeta$ is the coordinate in thickness direction $\vert \zeta \vert \le \sfrac{t}{2}$.

\subsection{Kinematics}
\label{sec:kin}

\begin{figure}[h]\centering
	
	% file=#1, size=#2, figno=#3
			\ifdefined\bUseTikzExternalize
				\def\tkzscale{0.275}
				\centering
				\tikzsetnextfilename{Fig2}  
				\input{tikz/Fig2}
			\else
				\includegraphics[width=1.0\textwidth]{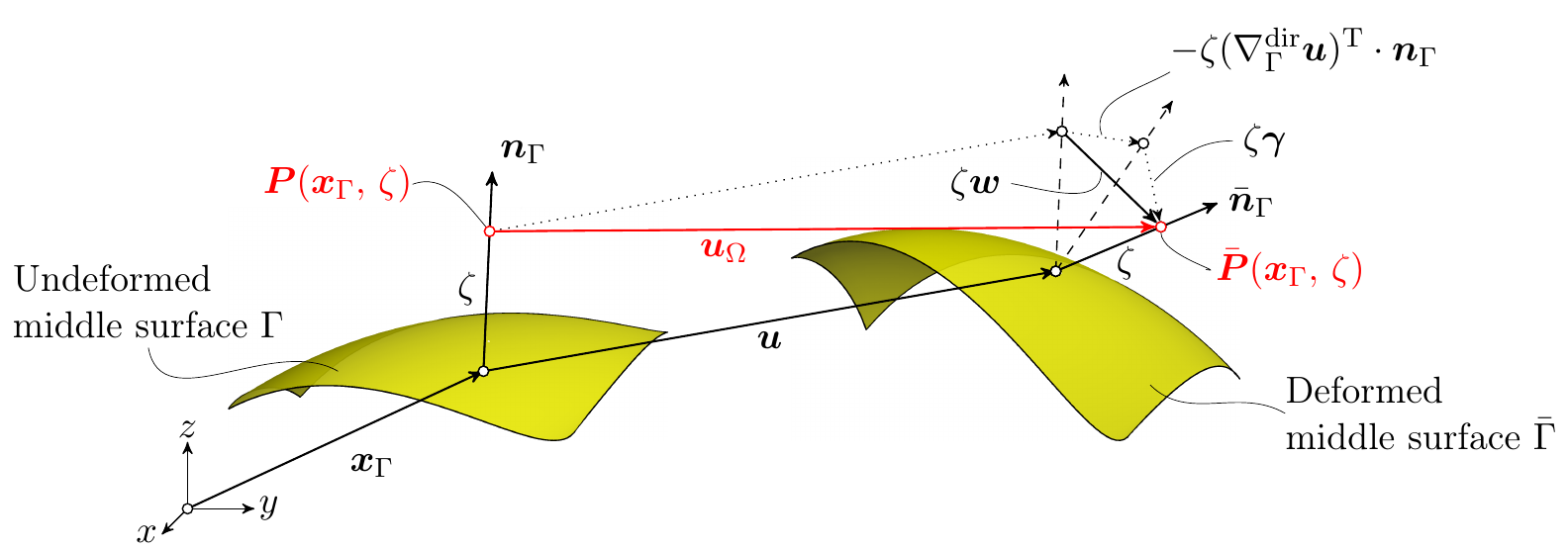}
			\fi			

	\caption{Displacement field $\vek{u}_\Omega$ of the Reissner-Mindlin shell.}
	\label{fig:dispu}
\end{figure}

For Reissner-Mindlin shells, the cross section remains straight after the deformation, but not necessarily normal to the middle surface due to transverse shear deformations. Herein, the rotation of the normal vector is modelled with a standard difference vector formulation \cite{Echter_2013a,Kiendl_2017a}. Other approaches such as exponential maps, rotation tensors, etc.~have been proposed, e.g., in \cite{Basar_1985a,Simo_1989a,Simo_1989b,Bischoff_2017a,Dornisch_2013a}, but are not considered here. The overall displacement of a point $\vek{P} \in \Omega$ is the difference between the spatial and reference configuration
\begin{align*}
\vek{u}_\Omega(\vek{x}) &= \bar{\vek{P}}(\vek{x}) - \vek{P}(\vek{x})\ ,
\end{align*}
which takes the form
\begin{align}
\vek{u}_\Omega(\vek{x}_\Gamma,\,\zeta) = \vek{u}(\vek{x}_\Gamma) + \zeta \vek{w}(\vek{x}_\Gamma)
\end{align}
with $\vek{u}(\vek{x}_\Gamma) : \Gamma \to \mathbb{R}^3$ being the displacement of the middle surface and $\vek{w}(\vek{x}_\Gamma)  : \Gamma \to T_P\Gamma$ being the difference vector, describing the rotation of the normal vector. In contrast to the Kirchhoff-Love shell, transverse shear deformations $\vek{\gamma}$ occur, which results in an additional rotation of the normal vector $\nG$, as illustrated in \autoref{fig:dispu}.
The difference vector $\vek{w}$ expressed in terms of the TDC is then defined as in \cite{Delfour_1996a,Schoellhammer_2018a} with additional transverse shear deformations
\begin{align}
\vek{w}(\vek{x}_\Gamma) = -[\gradGD{\vek{u}(\vek{x}_\Gamma)}]^\T \cdot \nG + \vek{\gamma}(\vek{x}_\Gamma)\ .
\end{align}
Note that the difference vector is a tangential vector as in the classical theory. The surface gradient of $\vek{u}_\Omega(\vek{x})$ is given by
\begin{align}
\begin{split}
\gradGD{\vek{u}_\Omega(\vek{x})} &= \mat{P} \cdot \dfrac{\p\vek{u}_\Omega(\vek{x})}{\p\vek{x}} \cdot \dfrac{\p \vek{x}}{\p \vek{x}_\Gamma}\ , \\
&= \mat{P}\cdot \left(\nabla\tilde{\vek{u}} + \nabla\zeta \otimes \vek{w} + \zeta \nabla\tilde{\vek{w}}\right) \cdot \left(\mathbb{I} + \zeta \mat{H}\right)\ ,\\
&=\left(\gradGD{\vek{u}} + \nG \otimes \vek{w} + \zeta \gradGD{\vek{w}}\right) \cdot \left(\mathbb{I} + \zeta \mat{H}\right)\ , \label{eq:graduom}
\end{split}
\end{align}
with $\mat{H} = \gradGD{\nG}$ being the Weingarten map \cite{Jankuhn_2017a,Delfour_2011a}. Note that the gradient of the thickness parameter $\zeta$ can be identified as the normal vector $\nG$. The second term $(\mathbb{I} + \zeta \mat{H})$ in \autoref{eq:graduom} is the inverse of the shell shifter \cite{Bischoff_2017a}, which is a second order tensor. \revStart The equivalent expression in local coordinates is $\sum_{i=1}^{3} \vek{G}^i \otimes \vek{A}_i$, where $\vek{G}^i$ are the contra-variant base vectors in the shell continuum and $\vek{A}_i$ the co-variant base vectors on the surface \cite{Bischoff_2017a}. \revEnd\par

The linearised strain tensor $\ten{\varepsilon}_\Gamma$ is defined by the symmetric part of the surface gradient of $\vek{u}_\Omega$
\begin{align}
\ten{\varepsilon}_\Gamma(\vek{x}) &= \dfrac{1}{2} \left[\gradGD{\vek{u}_\Omega} + (\gradGD{\vek{u}_\Omega})^\T\right] \\
&= \ten{\varepsilon}_\Gamma^\t{P}(\vek{x}) + \ten{\varepsilon}_\Gamma^\t{S}(\vek{x})
\end{align}
and is split into an in-plane strain $\ten{\varepsilon}_\Gamma^\t{P}$ and a transverse shear strain $\ten{\varepsilon}_\Gamma^\t{S}$. Neglecting higher-order terms in thickness direction, as usual in the classical theory \cite{Bischoff_2017a}, the in-plane strain is defined by \cite{Schoellhammer_2018a}
\begin{align}
\ten{\varepsilon}_\Gamma^\t{P} &= \mat{P} \cdot \ten{\varepsilon}_\Gamma \cdot \mat{P} \\
&= \ten{\varepsilon}_{\Gamma,\,\t{Mem}}^\t{P} + \zeta\ten{\varepsilon}_{\Gamma,\,\t{Bend}}^\t{P}\ ,
\end{align}
which is divided into a membrane and bending strain. The in-plane membrane strain becomes
\begin{align}
\ten{\varepsilon}_{\Gamma,\,\t{Mem}}^\t{P}(\vek{u}) = \dfrac{1}{2} \left[\gradGC{\vek{u}} + (\gradGC{\vek{u}})^\T \right]\ ,
\end{align}
and the bending strain is
\begin{align}
\ten{\varepsilon}_{\Gamma,\,\t{Bend}}^\t{P}(\vek{u}, \vek{w}) = \dfrac{1}{2} \left[ \mat{H}\cdot\gradGD{\vek{u}} + (\gradGD{\vek{u}})^\T\cdot \mat{H} + \gradGC{\vek{w}} + (\gradGC{\vek{w}})^\T\right]\ .
\end{align}
The transverse shear strain is defined in a similar manner as in \cite{Hansbo_2014b}
\begin{align}
\ten{\varepsilon}_\Gamma^\t{S}(\vek{u},\,\vek{w}) &= \mat{Q} \cdot \ten{\varepsilon}_\Gamma + \ten{\varepsilon}_\Gamma \cdot \mat{Q} \\
&= \dfrac{1}{2} \left[ \mat{Q} \cdot \gradGD{\vek{u}} + (\gradGD{\vek{u}})^\T \cdot \mat{Q} + \nG \otimes \vek{w} + \vek{w} \otimes \nG \right]\ .
\end{align}
When the shell surface is parametrized, the resulting strain components are equivalent compared to the classical theory in local coordinates, see, e.g., \cite{Basar_1985a,Bischoff_2017a,Echter_2013a,Kiendl_2017a}. In the case of flat shell structures, the membrane strain is only a function of the tangential part of the middle surface displacement $\vek{u}_t = \mat{P} \cdot\vek{u}$. Since the curvature is zero in this case, the Weingarten map $\mat{H}$ vanishes. Therefore, the bending strain is only a function of the difference vector $\vek{w}$ and the transverse shear strain becomes a function of the normal displacement $u_n = \vek{u} \cdot \nG$ and $\vek{w}$, resulting into the well-known Reissner-Mindlin \emph{plate}, see, e.g., \cite{Onate_2013a}.

\subsection{Constitutive Equation}
\label{sec:material}

The shell is assumed to be linear elastic and, as usual for thin structures, the Lam\'e constants are chosen such that the normal stress in thickness direction is eliminated, hence,
\begin{align}
\ten{\sigma}_\Gamma(\vek{x}) = 2\mu\ten{\varepsilon}_\Gamma(\vek{x}) + \lambda\t{tr}[\ten{\varepsilon}_\Gamma(\vek{x})]\mathbb{I}
\end{align}
where $\mu = \frac{E}{2(1+\nu)}$ and $\lambda = \frac{E\nu}{1-\nu^2}$. The stress tensor is decomposed in a similar manner as above into in-plane (membrane and bending)  stresses
\begin{align}\label{eq:inplanestress}
\ten{\sigma}_\Gamma^{\t{P}}(\vek{x}) = \mat{P} \cdot \left[2\mu\ten{\varepsilon}_\Gamma^{\t{P}}(\vek{x}) + \lambda\t{tr}[\ten{\varepsilon}_\Gamma^{\t{P}}(\vek{x})]\mathbb{I}\right] \cdot \mat{P}\ ,
\end{align}
and transverse shear stresses
\begin{align}
\begin{split}\label{eq:shearstress}
\ten{\sigma}_\Gamma^{\t{S}}(\vek{x}_\Gamma) &=  2\mu\left[\mat{Q} \cdot \ten{\varepsilon}_\Gamma^{\t{S}}(\vek{x}_\Gamma) + \ten{\varepsilon}_\Gamma^{\t{S}}(\vek{x}_\Gamma) \cdot \mat{Q}\right] + \lambda\t{tr}\left[\ten{\varepsilon}_\Gamma^{\t{S}}(\vek{x}_\Gamma) \right]\mat{Q}\\
&= 2\mu\,\alpha_\t{s}\, \ten{\varepsilon}_\Gamma^{\t{S}}(\vek{x}_\Gamma)\ .
\end{split}
\end{align}
As readily seen, the transverse shear stress is only a function of $\vek{x}_\Gamma$, which results in a constant transverse shear stress in thickness direction within the Reissner-Mindlin shell theory. In order account for this, a shear correction factor $\alpha_\t{s}$ is introduced \cite{Bischoff_2017a}. A common choice of the shear correction factor is $\alpha_\t{s} = \sfrac{5}{6}$.  Note that due to the double projection with $\mat{P}$ in \autoref{eq:inplanestress} of the in-plane stress, also directional gradients can be used, which is beneficial from an implementational point of view.

\subsection{Stress resultants}
\label{sec:stressres}

Assuming a constant shifter in the material law, the stress tensor $\ten{\sigma}_\Gamma(\vek{x})$ is only a function of the deflection of the middle surface $\vek{u}$, the difference vector $\vek{w}$ and linear in thickness direction. This enables an analytical pre-integration with respect to the thickness and stress resultants, such as effective membrane forces, bending moments and transverse shear forces can be identified. The following quantities are expressed in terms of the TDC using a global Cartesian coordinate system and are equivalent to the stress resultants in the classical theory using curvilinear coordinates, e.g., \cite{Basar_1985a,Bischoff_2017a}. \par

The symmetric, in-plane moment tensor $\mat{m}_\Gamma$ is defined as
\begin{align}
\mat{m}_\Gamma &= \int_{-\sfrac{t}{2}}^{\sfrac{t}{2}} \zeta \ten{\sigma}_\Gamma^\t{P}(\vek{x})\ \d\zeta = \dfrac{t^3}{12} \ten{\sigma}_\Gamma(\ten{\varepsilon}^\t{P}_{\Gamma,\t{Bend}}) \nonumber\\
&= \mat{P}\cdot {\mat{m}}^\t{dir}_\Gamma \cdot\mat{P}\ ,
\intertext{resulting in the components}
\left[{\mat{m}}^\t{dir}_\Gamma\right]_{11} &= D_{\t{B}}\, \left[w^\t{dir}_{x,x} + [\mat{H}]_{1j} \cdot \vek{u}^\t{dir}_{,x} + \nu(w^\t{dir}_{y,y} + w^\t{dir}_{z,z} + [\mat{H}]_{2j} \cdot \vek{u}^\t{dir}_{,y} + [\mat{H}]_{3j} \cdot \vek{u}^\t{dir}_{,z})\right] \nonumber\ ,\\[.1cm]
\left[{\mat{m}}^\t{dir}_\Gamma\right]_{22} &= D_{\t{B}}\, \left[w^\t{dir}_{y,y} + [\mat{H}]_{2j} \cdot \vek{u}^\t{dir}_{,y} + \nu(w^\t{dir}_{x,x} + w^\t{dir}_{z,z} + [\mat{H}]_{1j} \cdot \vek{u}^\t{dir}_{,x} + [\mat{H}]_{3j} \cdot \vek{u}^\t{dir}_{,z})\right]\nonumber\ ,\\[.1cm]
\left[{\mat{m}}^\t{dir}_\Gamma\right]_{33} &=D_{\t{B}}\, \left[w^\t{dir}_{z,z} + [\mat{H}]_{3j} \cdot \vek{u}^\t{dir}_{,z} + \nu(w^\t{dir}_{x,x} + w^\t{dir}_{y,y} + [\mat{H}]_{1j} \cdot \vek{u}^\t{dir}_{,x} + [\mat{H}]_{2j} \cdot \vek{u}^\t{dir}_{,y})\right]\nonumber\ ,\\[.1cm]
\left[{\mat{m}}^\t{dir}_\Gamma\right]_{12} &= D_{\t{B}}\, \frac{1-\nu}{2}\left[w^\t{dir}_{x,y} + w^\t{dir}_{y,x} + [\mat{H}]_{1j} \cdot \vek{u}^\t{dir}_{,y} + [\mat{H}]_{2j} \cdot \vek{u}^\t{dir}_{,x}\right]\nonumber\ ,\\[.1cm]
\left[{\mat{m}}^\t{dir}_\Gamma\right]_{13} &= D_{\t{B}}\, \frac{1-\nu}{2}\left[w^\t{dir}_{x,z} + w^\t{dir}_{z,x} + [\mat{H}]_{1j} \cdot \vek{u}^\t{dir}_{,z} + [\mat{H}]_{3j} \cdot \vek{u}^\t{dir}_{,x}\right]\nonumber\ ,\\[.1cm]
\left[{\mat{m}}^\t{dir}_\Gamma\right]_{23} &= D_{\t{B}}\, \frac{1-\nu}{2}\left[w^\t{dir}_{y,z} + w^\t{dir}_{z,y} + [\mat{H}]_{2j} \cdot \vek{u}^\t{dir}_{,z} + [\mat{H}]_{3j} \cdot \vek{u}^\t{dir}_{,y}\right]\nonumber\ ,
\end{align}
with $j = 1,\,2,\,3$ and $D_{\t{B}} = \frac{Et^3}{12(1-\nu^2)}$ being the flexural rigidity of the shell. The two non-zero eigenvalues of $\mat{m}_\Gamma$ are the principal moments $m_1$ and $m_2$. The effective membrane (normal) force tensor $\tilde{\ten{n}}_\Gamma$ is defined as
\begin{align}
\tilde{\mat{n}}_\Gamma &=\int_{-\sfrac{t}{2}}^{\sfrac{t}{2}}  \ten{\sigma}_\Gamma^\t{P}(\vek{x})\ \d\zeta = t \ten{\sigma}_\Gamma(\ten{\varepsilon}^\t{P}_{\Gamma,\t{Mem}}) \nonumber\\
&= \mat{P}\cdot {\mat{n}}^\t{dir}_\Gamma \cdot\mat{P}\ ,
\intertext{with components}
\left[{\mat{n}}^\t{dir}_\Gamma\right]_{11} &= D_{\t{M}} \left[u_{,x}^\t{dir} + \nu(v_{,y}^\t{dir}+w_{,z}^\t{dir})\right]\nonumber\ ,\\[.1cm]
\left[{\mat{n}}^\t{dir}_\Gamma\right]_{22} &= D_{\t{M}} \left[v_{,y}^\t{dir} + \nu(u_{,x}^\t{dir}+w_{,z}^\t{dir})\right]\nonumber\ ,\\[.1cm]
\left[{\mat{n}}^\t{dir}_\Gamma\right]_{33} &= D_{\t{M}} \left[w_{,z}^\t{dir} + \nu(u_{,x}^\t{dir}+v_{,y}^\t{dir})\right]\nonumber\ ,\\[.1cm]
\left[{\mat{n}}^\t{dir}_\Gamma\right]_{12} &=D_{\t{M}}\,\frac{1-\nu}{2}\left(u_{,y}^\t{dir}+v_{,x}^\t{dir}\right)\nonumber\ ,\\[.1cm]
\left[{\mat{n}}^\t{dir}_\Gamma\right]_{13} &= D_{\t{M}}\,\frac{1-\nu}{2}\left(u_{,z}^\t{dir}+w_{,x}^\t{dir}\right)\nonumber\ ,\\[.1cm]
\left[{\mat{n}}^\t{dir}_\Gamma\right]_{23} &= D_{\t{M}}\,\frac{1-\nu}{2}\left(v_{,z}^\t{dir}+w_{,y}^\t{dir}\right)\nonumber\ ,
\end{align}
where $D_{\t{M}} = \frac{Et}{1-\nu^2}$ is the membrane stiffness. Analogously to the moment tensor, the effective normal force tensor is also a symmetric, in-plane tensor. Note that for curved shells this tensor is \emph{not} the physical normal force tensor, but occurs in the weak form, see \autoref{sec:cwf}. In case of a curved shell, the physical normal force tensor $\mat{n}^{\t{real}}_\Gamma$ is defined by
\begin{align} \label{eq:nreal}
\mat{n}^{\t{real}}_\Gamma = \tilde{\mat{n}}_\Gamma + \mat{H}\cdot\mat{m}_\Gamma
\end{align}
and is, in general, not symmetric, but features one zero eigenvalue just as $\tilde{\mat{n}}_\Gamma$. With \autoref{eq:shearstress}, the resulting transverse shear force tensor is
\begin{align}
\mat{q}_\Gamma &=\int_{-\sfrac{t}{2}}^{\sfrac{t}{2}}  \ten{\sigma}_\Gamma^\t{S}(\vek{x})\ \d\zeta = t \ten{\sigma}_\Gamma(\ten{\varepsilon}^\t{S}_{\Gamma}) \nonumber\\
&=2D_\t{Shear}\,\ten{\varepsilon}_\Gamma^{\t{S}}\ ,
\intertext{with components}
\left[{\mat{q}}_\Gamma\right]_{11} &= 2D_{\t{Shear}} \left[n_x w_x + [\mat{Q}]_{1j} \cdot \vek{u}^\t{dir}_{,x} \right]\nonumber\ ,\\[.1cm]
\left[{\mat{q}}_\Gamma\right]_{22} &= 2D_{\t{Shear}} \left[n_y w_y + [\mat{Q}]_{2j} \cdot \vek{u}^\t{dir}_{,y} \right]\nonumber\ ,\\[.1cm]
\left[{\mat{q}}_\Gamma\right]_{33} &= 2D_{\t{Shear}} \left[n_z w_z + [\mat{Q}]_{3j} \cdot \vek{u}^\t{dir}_{,z} \right]\nonumber\ ,\\[.1cm]
\left[{\mat{q}}_\Gamma\right]_{12} &= D_{\t{Shear}} \left[n_x w_y + n_y w_x + [\mat{Q}]_{1j} \cdot \vek{u}^\t{dir}_{,y} + [\mat{Q}]_{2j} \cdot \vek{u}^\t{dir}_{,x} \right]\nonumber\ ,\\[.1cm]
\left[{\mat{q}}_\Gamma\right]_{13} &= D_{\t{Shear}} \left[n_x w_z + n_z w_x + [\mat{Q}]_{1j} \cdot \vek{u}^\t{dir}_{,z} + [\mat{Q}]_{3j} \cdot \vek{u}^\t{dir}_{,x} \right]\nonumber\ ,\\[.1cm]
\left[{\mat{q}}_\Gamma\right]_{23} &= D_{\t{Shear}} \left[n_y w_z + n_z w_y + [\mat{Q}]_{2j} \cdot \vek{u}^\t{dir}_{,z} + [\mat{Q}]_{3j} \cdot \vek{u}^\t{dir}_{,y} \right]\nonumber\ ,
\end{align}
where $D_{\t{Shear}} = \alpha_\t{s}\mu t = \alpha_\t{s}\frac{Et}{2(1+\nu)}$ is the transverse shear stiffness.

\subsection{Equilibrium in strong form}
\label{sec:equ}

Based on the stress resultants from above, one obtains the force and moment equilibrium for a curved Reissner-Mindlin shell in terms of the TDC using a global Cartesian coordinate system in strong form as
\begin{align}
\divG{\mat{n}^{\t{real}}_\Gamma} + \mat{Q}\cdot\divG{\mat{q}_\Gamma} + \mat{H}\cdot(\mat{q}_\Gamma\cdot\nG) &= - \vek{f}\ , \label{eq:sff}\\
\mat{P}\cdot\divG{\mat{m}_ {\Gamma}} - \mat{q}_\Gamma \cdot\nG &= -\vek{c}\ , \label{eq:sfm}
\end{align}
with $\vek{f}$ being the load vector per area and $\vek{c}$ being a distributed moment vector on the middle surface $\Gamma$. One may split the force equilibrium into the tangential and normal direction
\begin{align}
\mat{P}\cdot\divG{\mat{n}^{\t{real}}_\Gamma} + \mat{H}\cdot(\mat{q}_\Gamma\cdot\nG) &= - \vek{f}_t\ , \\
-\mat{H} : \mat{n}^{\t{real}}_\Gamma + \nG\cdot\divG{\mat{q}_\Gamma} &= - f_n\ .
\end{align}
Alternatively, \autoref{eq:sff} can be rewritten in terms of the effective normal force tensor by substituting $\mat{n}^{\t{real}}_\Gamma$ with \autoref{eq:nreal}
\begin{align}
\divG{\tilde{\mat{n}}_\Gamma} + \mat{H}\cdot\divG{\mat{m}_\Gamma} + \sum_{i,j=1}^{3} [\mat{H}_{,i}]_{jk} [\mat{m}_ {\Gamma}]_{ji} + \mat{Q}\cdot\divG{\mat{q}_\Gamma} + \mat{H}\cdot(\mat{q}_\Gamma\cdot\nG) &= - \vek{f}\ .
\end{align}

Assuming a bounded shell with boundary $\p\Gamma$, there exist for each field (deflection of the middle surface $\vek{u}$ and difference vector $\vek{w}$) two non-overlapping parts, the Dirichlet boundary $\p\Gamma_\t{D,i}$ and the Neumann boundary $\p\Gamma_\t{N,i}$, with $i = \lbrace\vek{u},\,\vek{w}\rbrace$. The corresponding boundary conditions for the displacement $\vek{u}$ are
\begin{align}
\begin{alignedat}{2}
\vek{u} &= \hat{\vek{g}}_{\vek{u}} &&\t{ on } \p\Gamma_{\t{D},\vek{u}}\ , \\
\mat{n}^{\t{real}}_\Gamma \cdot \nCo + (\nG\cdot\mat{q}_\Gamma\cdot\nCo)\cdot\nG &= \hat{\vek{p}} &&\t{ on } \p\Gamma_{\t{N},\vek{u}}\ .
\end{alignedat}
\end{align}
For the rotation of the normal vector, the boundary conditions are
\begin{align}
\begin{alignedat}{2}
\vek{w} &= \hat{\vek{g}}_{\vek{w}} &&\t{ on } \p\Gamma_{\t{D},\vek{w}}\ , \\
\mat{m}_\Gamma \cdot \nCo &= \hat{\vek{m}}_{\p\Gamma} &&\t{ on } \p\Gamma_{\t{N},\vek{w}}\ .
\end{alignedat}
\end{align}

\begin{figure}[h]\centering
	\subfloat[displacements and rotations]{
	% file=#1, size=#2, figno=#3
			\ifdefined\bUseTikzExternalize
				\def\tkzscale{0.2}
				\centering
				\tikzsetnextfilename{Fig3a}  
				\input{tikz/Fig3a}
			\else
				\includegraphics[width=0.45\textwidth]{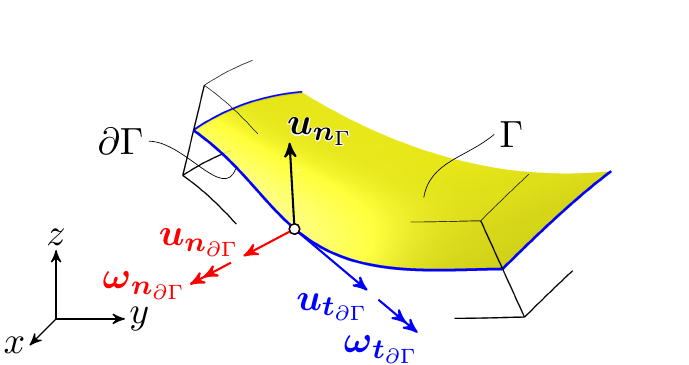}
			\fi			
\label{fig:bcrot}}
	\hfil
	\subfloat[forces and moments]{
	% file=#1, size=#2, figno=#3
			\ifdefined\bUseTikzExternalize
				\def\tkzscale{0.2}
				\centering
				\tikzsetnextfilename{Fig3b}  
				\input{tikz/Fig3b}
			\else
				\includegraphics[width=0.45\textwidth]{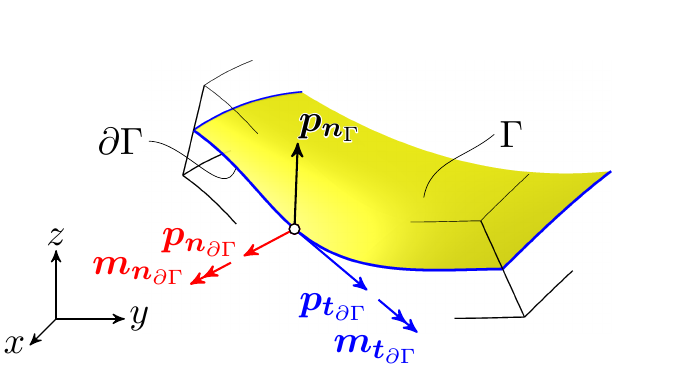}
			\fi			
\label{fig:bcnm}}	
	\caption{Decomposition of the middle surface displacement $\vek{u}$, difference vector $\vek{w}$, forces and bending moments along the boundary $\p\Gamma$ in terms of $\tB,\,\nCo$ and $\nG$: (a) displacements and rotations at the boundary, (b) forces and bending moments at the boundary.}	
	\label{fig:bcdecomp}
\end{figure}

In \autoref{fig:bcdecomp}, the possible boundary conditions are illustrated. In \autoref{fig:bcrot}, the displacement field $\vek{u}$ at the boundary is expressed in terms of the local triad $(\tB,\,\nCo,\,\nG)$ and, since the difference vector is tangential, the rotation of the normal vector may be written in terms of $(\tB,\,\nCo)$
\begin{align}
\vek{w} = \underset{\omega_{\tB}}{\underbrace{(\vek{w}\cdot\nCo)}}\, \tB + \underset{\omega_{\nCo}}{\underbrace{(\vek{w}\cdot\tB)}}\, \nCo\ .
\end{align}
The conjugated forces $(\vek{p}_{\tB},\,\vek{p}_{\nCo},\,\vek{p}_{\nG})$ and bending moments $(\vek{m}_{\tB},\,\vek{m}_{\nCo})$ at the boundary are shown in \autoref{fig:bcnm}.
A set of common support types is given in \autoref{tab:bc}. Other boundary conditions (e.g., membrane support, etc.) can be found, e.g., in \cite{Basar_1985a}.\par

\begin{table}[htb]
	\centering
	\begin{tabular}{lccccc}
		\toprule
		Clamped edge & ${u}_{\tB} = 0$ & ${u}_{\nCo} = 0$ & ${u}_{\nG} = 0$ & $\omega_{\tB} = 0$ &$\omega_{\nCo} = 0$\\ 
		Simply supported edge &${u}_{\tB} = 0$ & ${u}_{\nCo} = 0$ & ${u}_{\nG} = 0$ & $m_{\tB} = 0$ & $m_{\nCo} = 0$  \\ 		
		Symmetry support & ${p}_{\tB} = 0$ & ${u}_{\nCo} = 0$ & ${p}_{\nG} = 0$ & $\omega_{\tB} = 0$ & $m_{\nCo} = 0$ \\ 
		Free edge & ${p}_{\tB} = 0$ & ${p}_{\nCo} = 0$ & ${p}_{\nG} = 0$ & $m_{\tB} = 0$	& $m_{\nCo} = 0$ \\ \bottomrule		
	\end{tabular}
	\caption{Set of common boundary conditions}
	\label{tab:bc}
\end{table}
 
With the boundary conditions for the displacements and rotations, the complete second-order boundary value problem (BVP) is defined. The obtained BVP in terms of the TDC is valid in the case of implicitly and explicitly defined surfaces. In the case of parametrized shells, the equilibrium in strong form is equivalent to the strong form formulated in local coordinates \cite{Basar_1985a,Wempner_2002a,Kiendl_2017a}. However, because the obtained BVP does not rely on a parametrized middle surface of the shell, the formulation in the frame of TDC is more general.

\subsection{Equilibrium in weak form}
\label{sec:cwf}
We first formulate the weak form of the \emph{force} equilibrium. \autoref{eq:sff} is multiplied with a suitable test function $\vek{v}_{\vek{u}}$ and the divergence theorem \autoref{eq:DivTheorem} is applied to the terms $-\int_\Gamma \vu \cdot \divG{\mat{n}^{\t{real}}_\Gamma}\ \d A$ and $-\int_\Gamma \vu \cdot \mat{Q}\cdot\divG{\mat{q}_\Gamma}\ \d A$ resulting into
\begin{align}
\begin{split}\label{eq:wff}
\int_{\Gamma} \gradGD{\vu} : \tilde{\mat{n}}_\Gamma + (\mat{H}\cdot\gradGD{\vu}):\mat{m}_\Gamma + (\mat{Q}\cdot\gradGD{\vu}):\mat{q}_\Gamma\ \d A &= \\
&\hspace{-2.5cm} \int_{\Gamma} \vu \cdot \vek{f}\ \d A + \int_{\p\Gamma_{\t{N},\vek{u}}} \vu \cdot \hat{\vek{p}}\ \d s\
\end{split}
\end{align}
Note that, in order to obtain \autoref{eq:wff}  the identities $(\mat{H}\cdot\gradGD{\vu}):\mat{m}_\Gamma = \gradGD{\vu}:(\mat{H}\cdot\mat{m}_\Gamma)$ and $[\mat{Q}_{,x}^\Gamma\cdot\vu \quad \mat{Q}_{,y}^\Gamma\cdot\vu \quad \mat{Q}_{,z}^\Gamma\cdot\vu] : \mat{q}_\Gamma = \vu^\T \cdot \mat{H} \cdot (\mat{q}_\Gamma \cdot \nG)$ are used. As previously mentioned, in the weak form of the force equilibrium only the effective normal force tensor appears instead of the non-symmetric, physical normal force tensor.\par
The weak form of the \emph{moment} equilibrium is obtained in a similar manner by multiplying \autoref{eq:sfm} with a suitable, tangential test function $\vw$ and the divergence theorem of \autoref{eq:DivTheorem} is applied to the first term, leading to
\begin{align}
\begin{split}
\int_\Gamma \gradGD{\vw}:\mat{m}_\Gamma + \vw \cdot (\mat{q}_\Gamma\cdot\nG)\ \d \Gamma &=  \int_{\Gamma} \vw\cdot\vek{c}\ \d\Gamma + \int_{\p\Gamma_{\t{N},\vek{w}}} \vw \cdot \hat{\vek{m}}_{\p\Gamma}\ \d s\ .
\end{split}
\end{align}
Suitable trial and the test function spaces are subspaces of the $\mathcal{H}^1(\Gamma)^3$ Sobolev space, where $\mathcal{H}^1$ is the space of functions with square integrable first derivatives. 
\section{Implementational aspects}
\label{sec:impl}

The previously derived continuous weak forms can be discretized with different finite element approaches such as the classical Surface FEM or more recent approaches such as the CutFEM \cite{Burman_2015a,Burman_2018a,Cenanovic_2016a,Elfverson_2018a} and TraceFEM \cite{Grande_2016a,Olshanskii_2017a,Olshanskii_2017b,Reusken_2014a}. Herein, the weak form of the BVP is discretized using isogeometric analysis (IGA) as proposed by Hughes et al.~\cite{Hughes_2005a, Cottrell_2009a}, being closely related to Surface FEM when applied to solve PDEs on manifolds. It is pointed out that continuity requirements would also allow a standard FEM implementation using $C^0$-continuous shape functions. Nevertheless, we prefer the use of NURBS here, for example, due to the improved convergence properties and higher smoothness of the results (including smooth forces and moments). Furthermore, the improved smoothness enables us to compute errors in the strong form of the BVP, see the numerical results in \autoref{sec:numres}.\par

The NURBS patch $T$ is the middle surface of the shell and the elements $\tau_i\ ( i=1,\,\ldots,\,n_\t{Elem})$ are defined by the knot spans of the patch. Linking isogeometric analysis to standard FE terminology, one may naturally refer to (i) the NURBS patch as the ``mesh'', (ii) the knot spans as the ``elements'', and (iii) the NURBS-functions as the shape, test, and/or trial functions. The shape functions $N_{j}^k(\vek{x}_\Gamma(\vek{r}))$ employed for all (physical) fields involved are NURBS of order $k$ with $j = 1,\,\ldots,\,n_c$ being the number of control points. We avoid the mathematical definition of NURBS and the resulting patches here because of the abundance of literature devoted to IGA and consider this as common state of the art.\par

The surface derivatives of the shape functions $\gradG{N}(\vek{x}_\Gamma(\vek{r}))$ are computed as defined in \autoref{eq:TangGradExplicit} , similar to the Surface FEM \cite{Dziuk_2013a,Demlow_2009a,Fries_2018b,Fries_2017b} using NURBS instead of Lagrange polynomials as trial and test functions. A general finite element space of order $k$ is now defined by
\begin{align*}
\mathcal{Q}^h_k = \left\lbrace u^h \in C_{k-1}(\Gamma), u^h = \sum_{j = 1}^{n_c} N_{j}^k(\vek{x}_\Gamma(\vek{r})) \hat{u}_j\,,\ \hat{u}_j \in \mathbb{R}\right\rbrace \subset \mathcal{H}^1(\Gamma)\ ,
\end{align*}
where the degrees of freedom $\hat{u}_j$ are stored at the control points.\par

The discrete displacement of the middle surface results as $ \vek{u}^h =  u^{h,i}\vek{E}_i$, with $\vek{E}_i$ being Cartesian base vectors, with $i = 1,\,2,\,3$ and $u^{h,i} = \Nu^\T \cdot \hat{\vek{u}}^i$. In contrast to $\vek{u}$, the difference vector $\vek{w}$ is a tangential vector, describing the rotation of the normal vector. The discretization of a tangential vector is, in general, not straightforward and, in the following, different strategies are examined: \par

(1) In the case of a parametrized surface, the co-variant base vectors $\vek{A}_\alpha,\ \alpha =1,\,2$, which are by construction tangential, may be used to define the difference vector $\vek{w}^h = w^{h,\,\alpha} \vek{A}_\alpha$, 
where $w^{h,\,\alpha} = \Nw^\T \cdot \hat{\vek{w}}^\alpha$. This approach is used in the classical 5-parameter models \cite{Bischoff_2017a}.\par

(2) Alternatively, the directions of the principal curvatures, which are eigenvectors of the Weingarten map $\mat{H}$, can be used as basis vectors. These vectors are perpendicular and also tangential by construction. This might be a reasonable choice in the case of curved, implicitly defined surfaces, where a parametrization is not available. Compared to the first approach, the crucial requirement of a parametrization is circumvented and the number of degrees of freedom per control point is equal.\par

(3) Another possibility is to define the difference vector in the global Cartesian coordinate system $\vek{w}^h = w^{h,i}\vek{E}_i$, with $i = 1,\,2,\,3$ and $w^{h,i} = \Nw^\T \cdot \hat{\vek{w}}^i$ and the constraint $\vek{w}^h \cdot \nG = 0$ is weakly enforced with a Lagrange multiplier or with the penalty method.\par

(4) A variant of (3), is to project the difference vector onto the tangent space of the middle surface $\vek{w}^h = \mat{P} \cdot w^{h,i}\vek{E}_i$. An advantage of this approach is, that the additional Lagrange multiplier field is not needed. On the other hand, due to the projection, conditioning issues occur, which may be addressed with an additional stabilization term.\par

Herein the third approach, where the difference vector is globally defined and the constraint is enforced with a Lagrange multiplier, is chosen. The shape functions of the discrete Lagrange multiplier $\lambda^h_n = \vek{N}^\T_{\lambda_n} \cdot \hat{\vek{\lambda}}_n$ for the constraint on the difference vector is defined in the same  manner as the components of the middle surface displacement. Furthermore, the boundary conditions shall be enforced weakly with Lagrange multipliers \cite{Zienkiewicz_2013a}. The shape functions of the discrete Lagrange multiplier field for the displacements is defined as ${N}_{\lambda_{\vek{u}}} = \lbrace {N}_{{\vek{u}}} \rvert_{\p\Gamma_{\t{D},\vek{u}}} \rbrace$ and for the difference vector as  ${N}_{\lambda_{\vek{w}}} = \lbrace {N}_{{\vek{w}}} \rvert_{\p\Gamma_{\t{D},\vek{w}}} \rbrace$.\par

Based on this, the following discrete trial and test functions spaces are defined:
\begin{align} \label{eq:funb}
\mathcal{S}_{\vek{u}}^h &= \mathcal{V}_{\vek{u}}^h = \mathcal{S}_{\vek{w}}^h = \mathcal{V}_{\vek{w}}^h = \left\lbrace \vek{u}^h \in \left[\mathcal{Q}^h_k\right]^3\right\rbrace \\
\mathcal{L}_{\lambda_n}^h &= \mathcal{V}_{\lambda_n}^h = \left\lbrace {\lambda}_n^h \in \left[\mathcal{Q}^h_k\right]\right\rbrace \label{eq:funlm}\\
\mathcal{L}_{\vek{\lambda}_{\vek{u}}}^h &= \mathcal{V}_{\vek{\lambda}_{\vek{u}}}^h = \left\lbrace \vek{\lambda}_{\vek{u}}^h\rvert_{\p\Gamma_{\t{D},\vek{u}}} : \ \vek{\lambda}_{\vek{u}}^h \in \mathcal{S}_{\vek{u}}^h \right\rbrace \\
\mathcal{L}_{\vek{\lambda}_{\vek{w}}}^h &= \mathcal{V}_{\vek{\lambda}_{\vek{w}}}^h = \left\lbrace \vek{\lambda}_{\vek{w}}^h\rvert_{\p\Gamma_{\t{D},\vek{w}}} : \ \vek{\lambda}_{\vek{w}}^h \in \mathcal{S}_{\vek{w}}^h \right\rbrace \label{eq:fune}
\end{align}

\subsection{Discrete weak form}
\label{sec:dwf}

The discrete weak form of the Reissner-Mindlin shell with Lagrange multipliers is the following. Given Young's modulus $E \in \mathbb{R}^+$, Poisson's ratio $\nu \in [0,\,0.5]$, surface load and moment $\vek{f}^h,\,\vek{c}^h$ on $\Gamma$, traction $\hat{\vek{p}}^h$ on
$\p\Gamma_{\t{N},\vek{u}}$, bending moments $\hat{\vek{m}}_{\p\Gamma}^h$ on $\p\Gamma_{\t{N},\vek{w}}$ and boundary conditions $\hat{\vek{g}}_{\vek{u}}^h$ in $\p\Gamma_{\t{D},\vek{u}}$, $\hat{\vek{g}}_{\vek{w}}^h$ on $\p\Gamma_{\t{D},\vek{w}}$, find the displacement field $\vek{u}^h \in \mathcal{S}_{\vek{u}}^h$, the difference vector $\vek{w}^h \in \mathcal{S}_{\vek{w}}^h$, and the Lagrange multiplier fields $(\lambda^h_n,\,\vek{\lambda}_{\vek{u}}^h,\,\vek{\lambda}_{\vek{w}}^h) \in \mathcal{L}_{\lambda_n}^h \times \mathcal{L}_{\vek{\lambda}_{\vek{u}}}^h \times \mathcal{L}_{\vek{\lambda}_{\vek{w}}}^h$ such that for all test functions $(\vu^h,\,\vw^h,\,v^h_{\lambda_n},\,\vek{v}^h_{\vek{\lambda}_{\vek{u}}} ,\,\vek{v}^h_{\vek{\lambda}_{\vek{w}}} ) \in \mathcal{V}_{\vek{u}}^h \times \mathcal{V}_{\vek{w}}^h \times \mathcal{V}_{\lambda_n}^h \times \mathcal{V}_{\vek{\lambda}_{\vek{u}}}^h \times \mathcal{V}_{\vek{\lambda}_{\vek{w}}}^h$, there holds in $\Gamma$
\begin{align*}
\int\limits_{\Gamma} \gradGD{\vu^h} : \tilde{\mat{n}}_\Gamma + (\mat{H}\cdot\gradGD{\vu^h}):\mat{m}_\Gamma + (\mat{Q}\cdot\gradGD{\vu^h}):\mat{q}_\Gamma\ \d A &= \int\limits_{\Gamma} \vu^h \cdot \vek{f}^h\ \d A + \!\!\!\int\limits_{\p\Gamma_{\t{N},\vek{u}}}  \!\!\!\! \vu^h \cdot \hat{\vek{p}}^h\ \d s\,,\\
\int\limits_\Gamma \gradGD{\vw^h}:\mat{m}_\Gamma + \vw^h \cdot (\mat{q}_\Gamma\cdot\nG) + \lambda^h_n(\vw^h \cdot \nG) \ \d A&= \int\limits_{\Gamma} \vw^h\cdot\vek{c}^h\ \d\Gamma + \!\!\!\!\int\limits_{\p\Gamma_{\t{N},\vek{w}}} \!\!\!\! \vw^h \cdot \hat{\vek{m}}^h_{\p\Gamma}\ \d s\,,\\
\int\limits_{\Gamma} v^h_{\lambda_n}\vek{w}^h \cdot \nG\ \d A &= 0\,,\\
\int\limits_{\p\Gamma_{\t{D},\vek{u}}} \vek{v}_{\vek{\lambda}_{\vek{u}}}^h \cdot \vek{u}^h\ \d s &= \int\limits_{\p\Gamma_{\t{D},\vek{u}}} \vek{v}_{\vek{\lambda}_{\vek{u}}}^h \cdot \vek{g}_{\vek{u}}^h\ \d s\,,\\
\int\limits_{\p\Gamma_{\t{D},\vek{w}}} \vek{v}_{\vek{\lambda}_{\vek{w}}}^h \cdot \vek{w}^h\ \d s &= \int\limits_{\p\Gamma_{\t{D},\vek{w}}} \vek{v}_{\vek{\lambda}_{\vek{w}}}^h \cdot \vek{g}_{\vek{w}}^h\ \d s\,.
\end{align*}
The usual element assembly yields a linear system of equations (if displacements and rotations are prescribed) in the form
\begin{align}\label{eq:linsys}
\begin{bmatrix}
\mat{K}_{\vek{u}\vek{u}} & \mat{K}_{\vek{u}\vek{w}} & \mat{0} & \mat{L}_{\lambda_{\vek{u}}} & \mat{0}\\
 \mat{K}^\T_{\vek{u}\vek{w}} &  \mat{K}_{\vek{w}\vek{w}} &  \mat{L}_{\lambda_n} & \mat{0} & \mat{L}_{\lambda_{\vek{w}}}\\ 
\mat{0} &\mat{L}_{\lambda_n}^\T & \mat{0} & \mat{0} & \mat{0}\\
\mat{L}_{\lambda_{\vek{u}}}^\T & \mat{0} & \mat{0} & \mat{0} & \mat{0} \\
\mat{0} & \mat{L}_{\lambda_{\vek{w}}}^\T & \mat{0} & \mat{0} & \mat{0}
\end{bmatrix} \cdot \begin{bmatrix}
\underline{\hat{\vek{u}}} \\
\underline{\hat{\vek{w}}} \\
\hat{\vek{\lambda}}_n \\
\underline{\hat{\vek{\lambda}}}_{\vek{u}} \\
\underline{\hat{\vek{\lambda}}}_{\vek{w}}
\end{bmatrix} = \begin{bmatrix}
\vek{b}_{\vek{f}} \\
\vek{b}_{\vek{c}} \\
\vek{0} \\
\vek{b}_{\vek{\lambda}_{\vek{u}}} \\
\vek{b}_{\vek{\lambda}_{\vek{w}}}
\end{bmatrix}\ ,
\end{align}
with $[\underline{\hat{\vek{u}}},\ \underline{\hat{\vek{w}}},\ \hat{\vek{\lambda}}_n,\ \underline{\hat{\vek{\lambda}}}_{\vek{u}},\ \underline{\hat{\vek{\lambda}}}_{\vek{w}}]^\T$ being the sought displacements and rotations of the control points and Lagrange multiplier fields. As usual in the context of Lagrange multiplier methods, \autoref{eq:linsys} has a saddle point structure and the well-known Babu{\v{s}}ka-Brezzi condition \cite{Babuska_1971a,Brezzi_1974,Franca_1988a} must be satisfied in order to obtain useful solutions in all involved fields. \par \revStart

We have conducted a number of numerical experiments with respect to the orders of the approximations of the displacements, rotations, and the Lagrange multiplier fields enforcing the boundary conditions and constraint on the difference vector. Results presented in the next section are obtained using the same order of these fields, which is implementationally the most convenient setting. Choosing the order for the constraint on the difference vector one order higher or lower did not noticeably change these results. We have also studied the penalty method for this constraint and obtained optimal convergence rates for the penalty parameter $\alpha$ in the range between $10^6$ and $10^{10}$.\par
\revEnd

\section{Numerical results}
\label{sec:numres}

In this section the proposed reformulation of the classical Reissner-Mindlin shell equations in the frame of the TDC is applied to a set of benchmark examples, consisting of the well-known Scordelis-Lo roof from \cite{Belytschko_1985a}, the partly clamped hyperbolic paraboloid from \cite{Bathe_2000a,Chapelle_1998a} and the clamped flower shaped shell from \cite{Schoellhammer_2018a}. In all examples, the shells are rather thin and locking phenomena can be expected, especially in the case of low ansatz orders. However, when increasing the order $p$, locking phenomena decrease significantly and, therefore, no further measures against locking phenomena are considered herein.\par

In the convergence studies, uniform NURBS patches with different orders and numbers of knot
spans in each direction are employed. This is equivalent to meshes with higher-order
elements and $n = \lbrace 2,\,4,\,8,\,16,\,32,\,64,\,128 \rbrace$ elements per side are used. The orders are varied as $p = \lbrace2,\,3,\,4,\,5,\,6\rbrace$.

\subsection{Scordelis-Lo roof}
\label{sec:scordelis}

\begin{figure}[ht]
	\begin{minipage}{.5\textwidth}\centering
		\includegraphics[width=.95\textwidth, height = .9\textwidth, keepaspectratio]{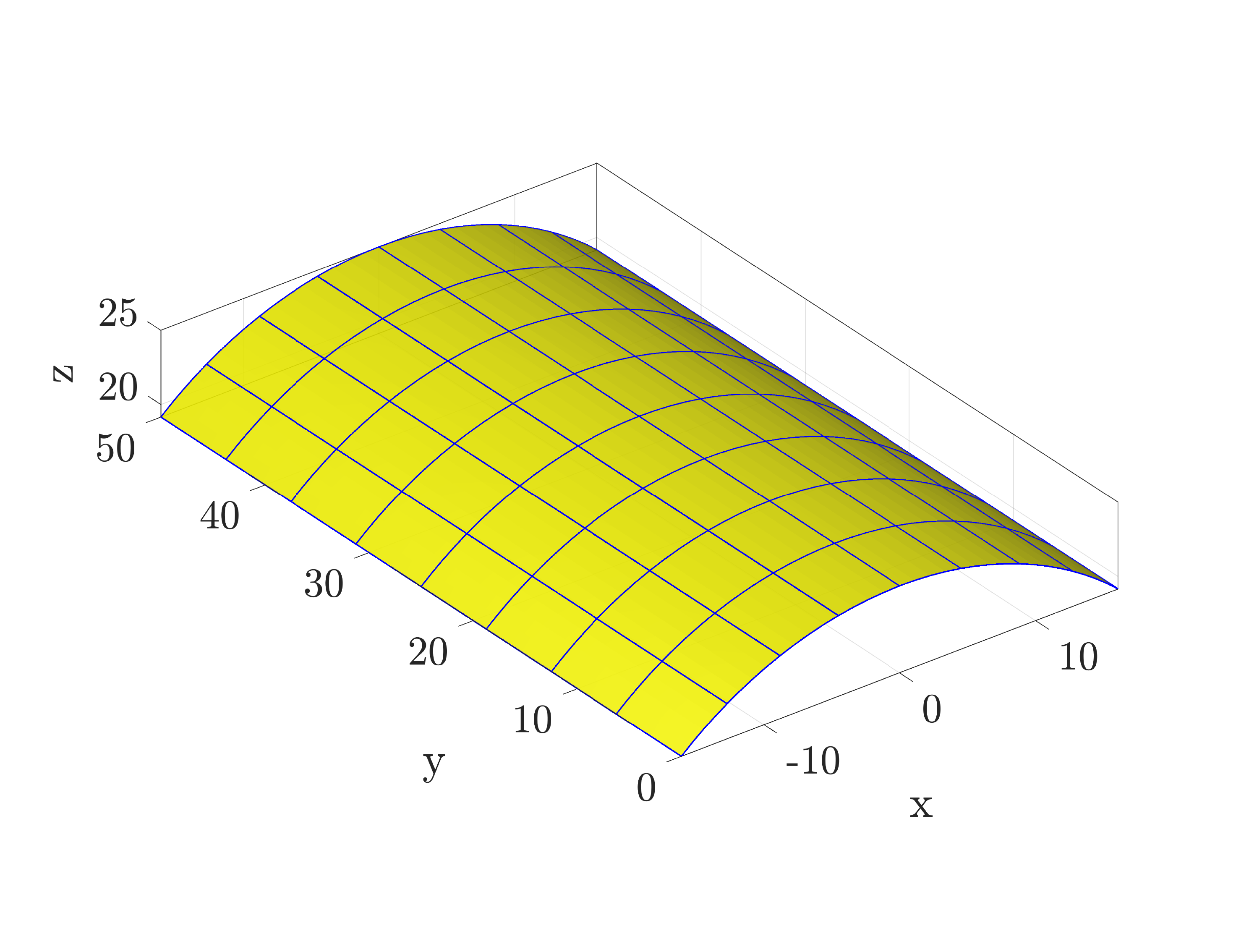}
	\end{minipage}
	\begin{minipage}{.48\textwidth}\scriptsize\flushleft
		\begin{tabular}[h]{ll}
			Geometry: & Cylindrical shell\\
			& $L = 50$ \\
			& $R = 25 $ \\
			& $\phi = \SI{80}{\degree}$ \\
			& $t = \SI{0.25}{}$\\[.25 cm]
			Material parameters: & $E = \SI{4.32e8}{}$ \\
			& $\nu = \SI{0.0}{}$\\
			& $\alpha_{\t{s}} = 1.0$ \\[.25 cm]
			Load: & Gravity load $\vek{f} = [0,\,0,\,-90]^\T $\\[.1cm]
			& \phantom{Gravity load} $\vek{c} = \vek{0}$\\[.25 cm]
			Support: & Rigid diaphragms at it ends
		\end{tabular}
	\end{minipage}
	\caption{Definition of Scordelis-Lo roof problem.}
	\label{fig:overscor}
\end{figure}

The Scordelis-Lo roof is part of the so-called \emph{shell obstacle course} \cite{Belytschko_1985a}. The shell is a cylindrical section and is supported with two rigid diaphragms at the ends and loaded by gravity forces, see \autoref{fig:overscor}. 
The cylinder is defined with a length $L = \SI{50}{}$, a radius of $R = \SI{25}{}$ and the angle subtended by the roof is $\phi = \SI{80}{\degree}$. The thickness of the shell is set to $t = \SI{0.25}{}$. The material parameters are Young's modulus $E = \SI{4.32e8}{}$ and the Poisson's ratio $\nu = \SI{0.0}{}$. In the convergence study, the maximum vertical displacement $u_{z,\max}$ is compared with the reference solution $u_{z,\max,\text{Ref}} = \SI{0.3024}{}$ as given in \cite{Belytschko_1985a}. The largest vertical displacement occurs in the midpoint of the free edges $[\pm R\cos(\SI{50}{\degree}),\,25,\, R\sin(\SI{50}{\degree})]^\T$.

\begin{figure}[ht]
	\centering
	\subfloat[displacement $\vek{u}$]{\includegraphics[width=.52\textwidth, height = .52\textwidth, keepaspectratio]{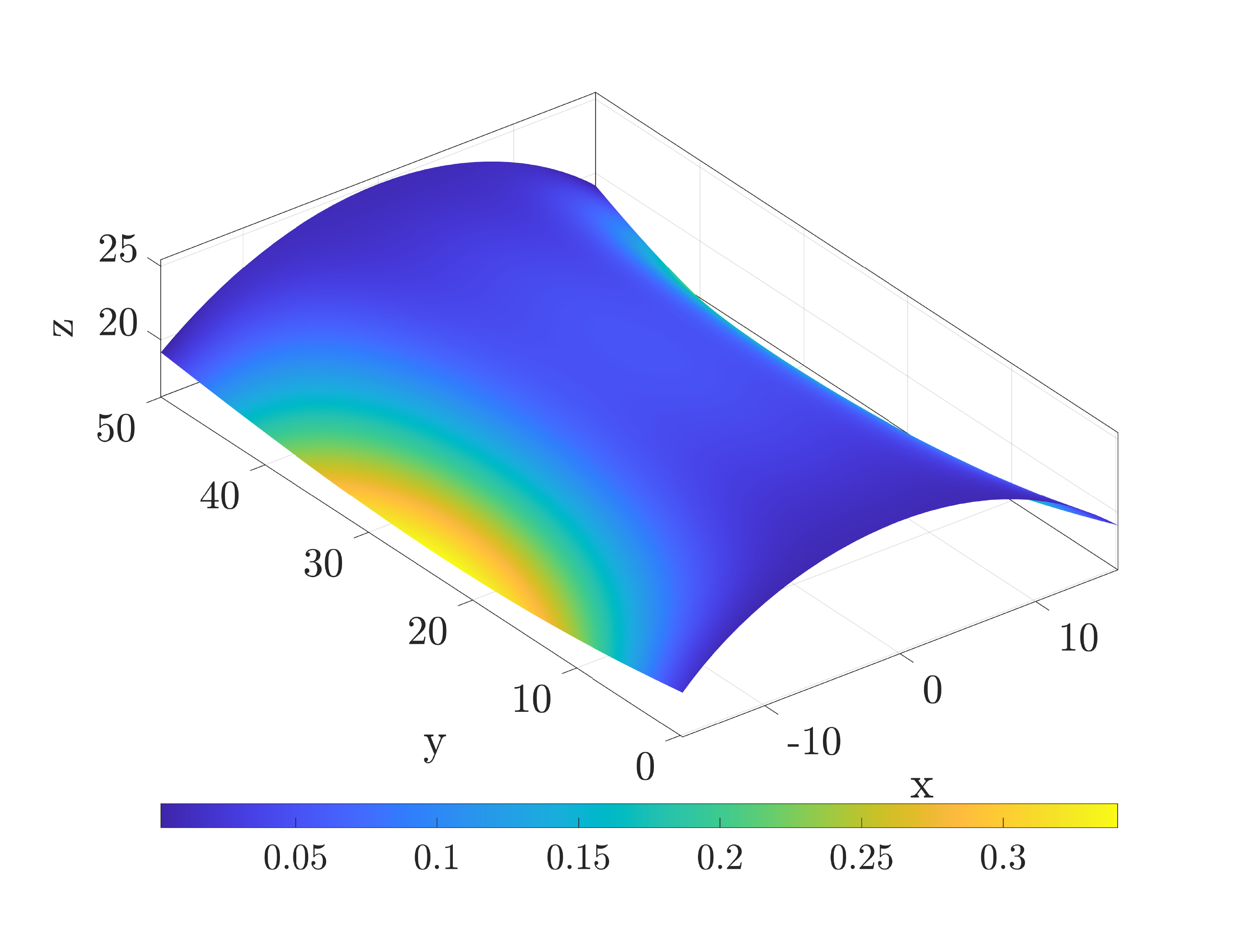}\label{fig:scordelisdisp}}
	\hfil
	\subfloat[convergence]{\includegraphics[width=.45\textwidth, height = .4\textwidth, keepaspectratio]{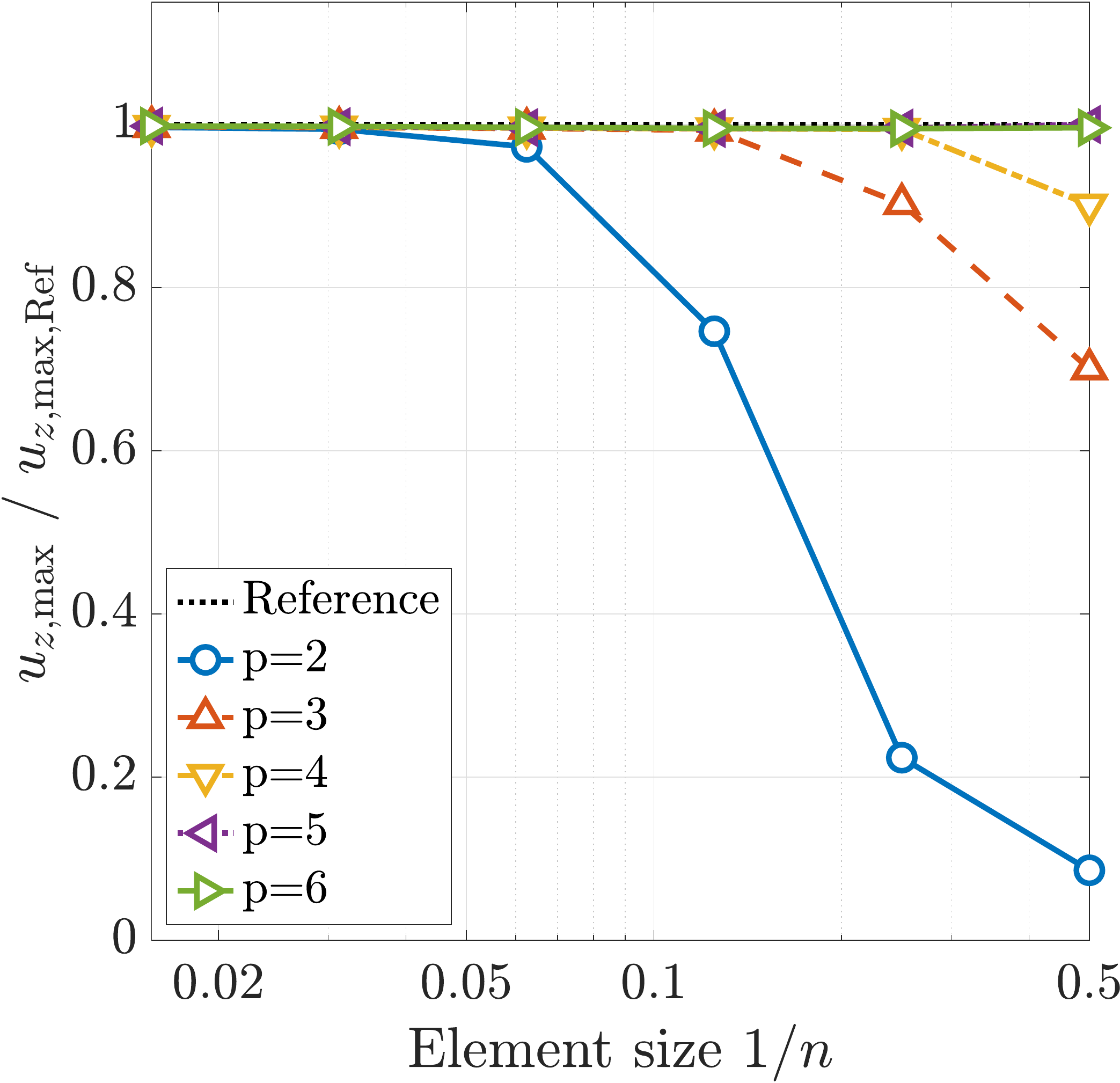}\label{fig:scordelisl2}}
	\caption{(a) Displacement field of the Lo-Scordelis roof scaled by one order of magnitude, (b) normalized convergence of reference displacement $u_{z,\max,\text{Ref}} = \SI{0.3024}{}$.}
\end{figure}

In \autoref{fig:scordelisdisp}, the numerical solution of the Scordelis-Lo roof is illustrated. The displacements are magnified by one order of magnitude. The colors on the deformed surface indicate the Euclidean norm of the displacement field $\vek{u}$.  In \autoref{fig:scordelisl2}, the normalized convergence of the maximum displacement $u_{z,\max}$ is plotted up to polynomial order of $p=6$ as a function of the element (knot span) size. \revStart It is clearly seen that the results improve upon increasing the order of the NURBS. It is noted, that due to the lack of smoothness of the solution for this classical benchmark test (at the corners of the shell), optimal higher-order convergence rates (as typically visualized in double-logarithmic error plots) can not be expected. Hence, the visualization as in \autoref{fig:scordelisl2} follows the usual style of many other references such as, e.g., in \cite{Belytschko_1985a,Cirak_2000a,Kiendl_2009a}. \revEnd Except for the order $p = 2$, the locking phenomena are not very pronounced. Hence, graded meshes as used, e.g., in \cite{Kiendl_2017a} in order to resolve the boundary layers, are not used here for the sake of simplicity. 

\subsection{Partly clamped hyperbolic paraboloid}
\label{sec:hyperbolic}

\begin{figure}[ht]
	\begin{minipage}{.5\textwidth}\centering
		\includegraphics[width=.95\textwidth, height = .9\textwidth, keepaspectratio]{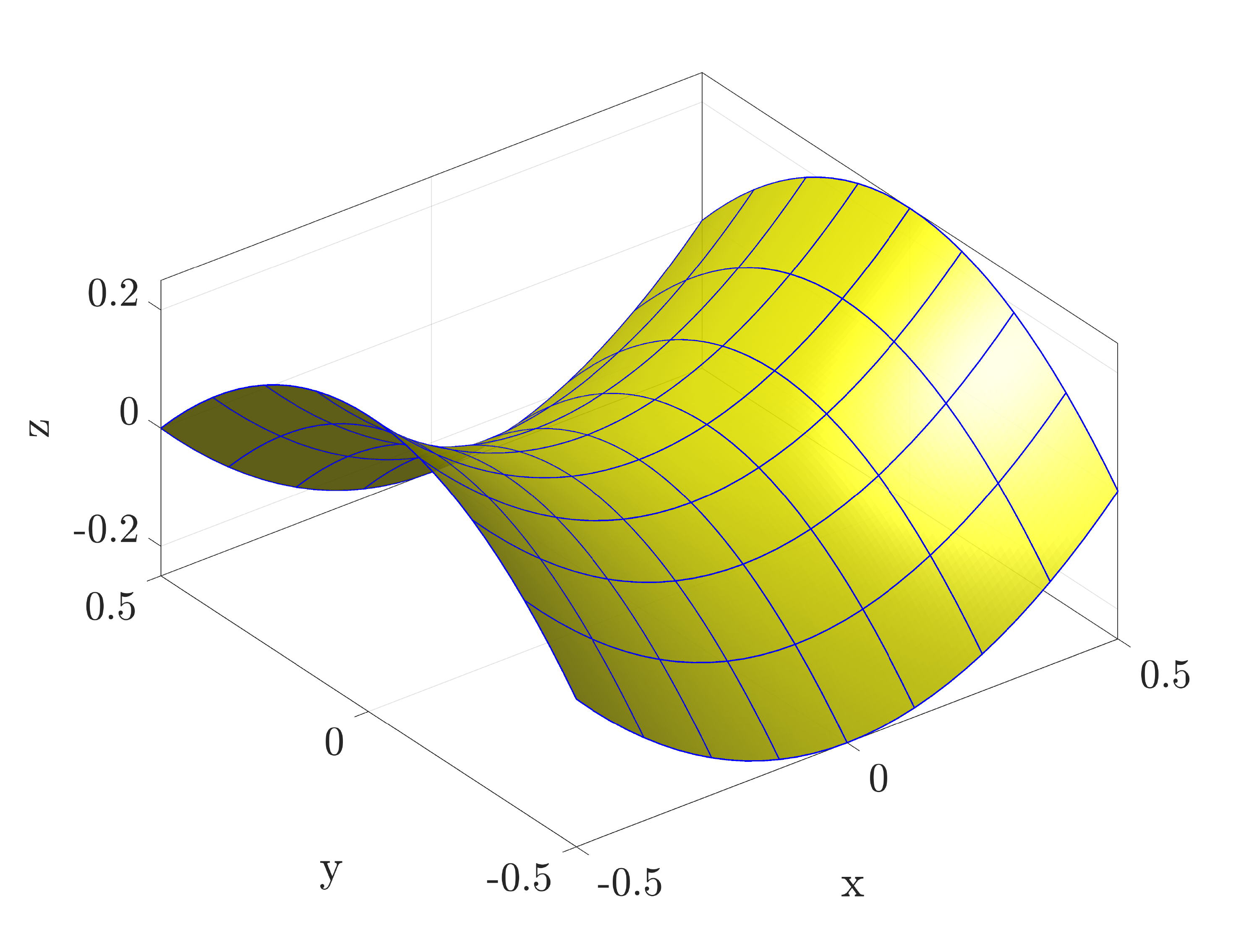}
	\end{minipage}
	\begin{minipage}{.48\textwidth}\scriptsize\flushleft
		\begin{tabular}[h]{ll}
			Geometry: & Hyperbolic paraboloid\\
			& $L_x = L_y = 1$ (center is at origin) \\
			& $z = x^2 - y^2 $ \\
			& $t = \SI{0.01}{}$\\[.25 cm]
			Material parameters: & $E = \SI{2.0e11}{}$ \\
			& $\nu = \SI{0.3}{}$\\
			& $\alpha_{\t{s}} = 1.0$ \\[.25 cm]
			Load: & Gravity load $\vek{f} = [0,\,0,\,-8000\cdot t]^\T $\\[.1 cm]
						& \phantom{Gravity load} $\vek{c} = \vek{0}$\\[.25 cm]
			Support: & Clamped edge at $x = -\sfrac{1}{2}$
		\end{tabular}
	\end{minipage}
	\caption{Definition of partly clamped hyperbolic paraboloid problem.}
	\label{fig:overhyper}
\end{figure}

The next test case is a partly clamped hyperbolic paraboloid and is taken from \cite{Bathe_2000a,Chapelle_1998a}. The shell is defined by $z = x^2 - y^2$ with $ (x,\,y) \in [-\sfrac{1}{2},\,\sfrac{1}{2}]^2$, the thickness is set to $t = 0.01$ and is loaded by gravity forces, see \autoref{fig:overhyper}. The edge at $x = - \sfrac{1}{2}$ is clamped and the other three edges are free. The material parameters are Young's modulus $E = \SI{2.0e11}{}$, Poisson's ratio $\nu = \SI{0.3}{}$. Similar to the example before, the displacements are compared with a reference solution. In particular, the vertical displacement at point $i = (0.5,\,0,\,0.25)^\T$ is compared with the reference solution $u_{z,\t{Ref}} = \SI{-9.3355e-5}{}$ given in \cite{Bathe_2000a}.

\begin{figure}[ht]
	\centering
	\subfloat[displacement $\vek{u}$]{\includegraphics[width=.52\textwidth, height = .52\textwidth, keepaspectratio]{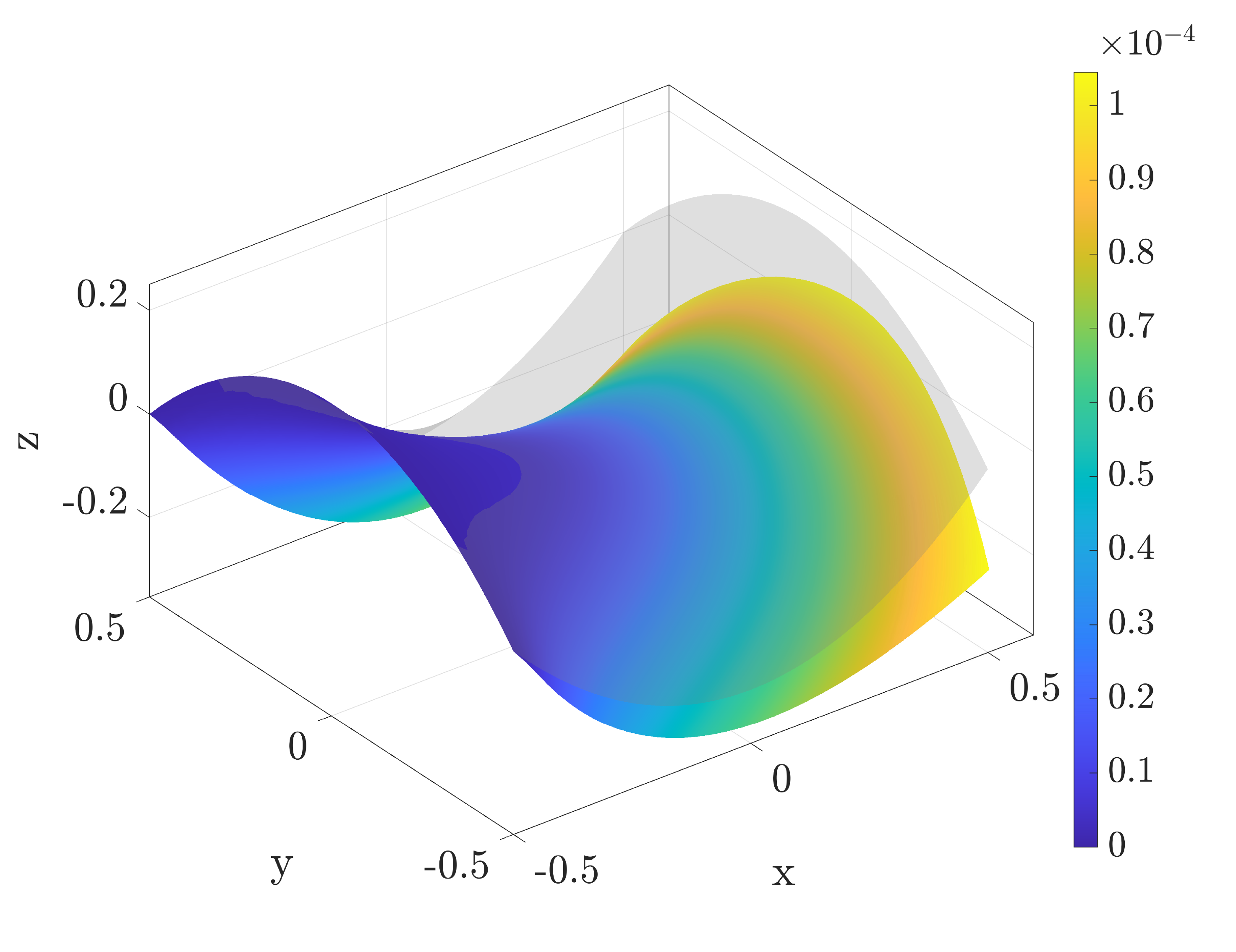}\label{fig:hyperdisp}}
	\hfil
	\subfloat[convergence]{\includegraphics[width=.45\textwidth, height = .45\textwidth, keepaspectratio]{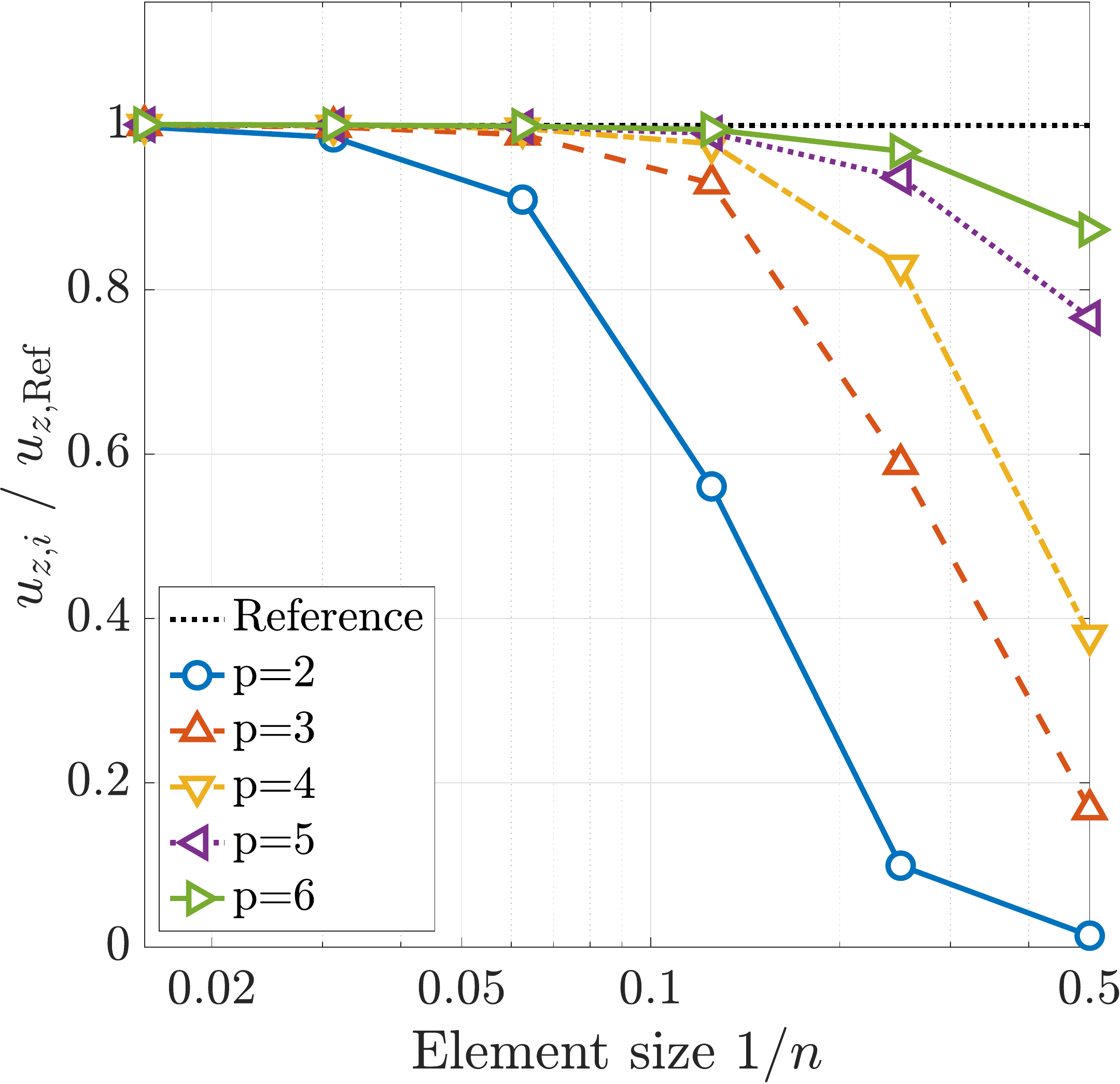}\label{fig:hyperl2}}
	\caption{(a) Displacement field of the partly clamped hyperbolic paraboloid scaled by  $\SI{2000}{}$, (b) normalized convergence of reference displacement $u_{z,\text{Ref}} = \SI{-9.3355e-5}{}$ at point $i = (0.5,\,0,\,0.25)^\T$.}
\end{figure}

In \autoref{fig:hyperdisp}, the undeformed domain (grey) and the deformed shell is presented, with displacements scaled by a factor of $\SI{2000}{}$. In \autoref{fig:hyperl2}, \revStart analogously to the example before, \revEnd the normalized convergence of the vertical displacement $u_{z,i}$ at point $i$ is plotted up to polynomial order of $p=6$ as a function of the element (knot span) size. For the lower orders $p=2,\,3$, the expected locking phenomena is more pronounced compared to the example before. Nevertheless, it is clearly seen that the accuracy for higher-order NURBS increases significantly and the behaviour of convergence is in agreement with the results shown e.g., in \cite{Bathe_2000a,Kiendl_2017a}. 

\subsection{Clamped flower shaped shell}
\label{sec:flower}

\begin{figure}[hbt]
	\begin{minipage}{.5\textwidth}\centering
		\includegraphics[width=.9\textwidth, height = .9\textwidth, keepaspectratio]{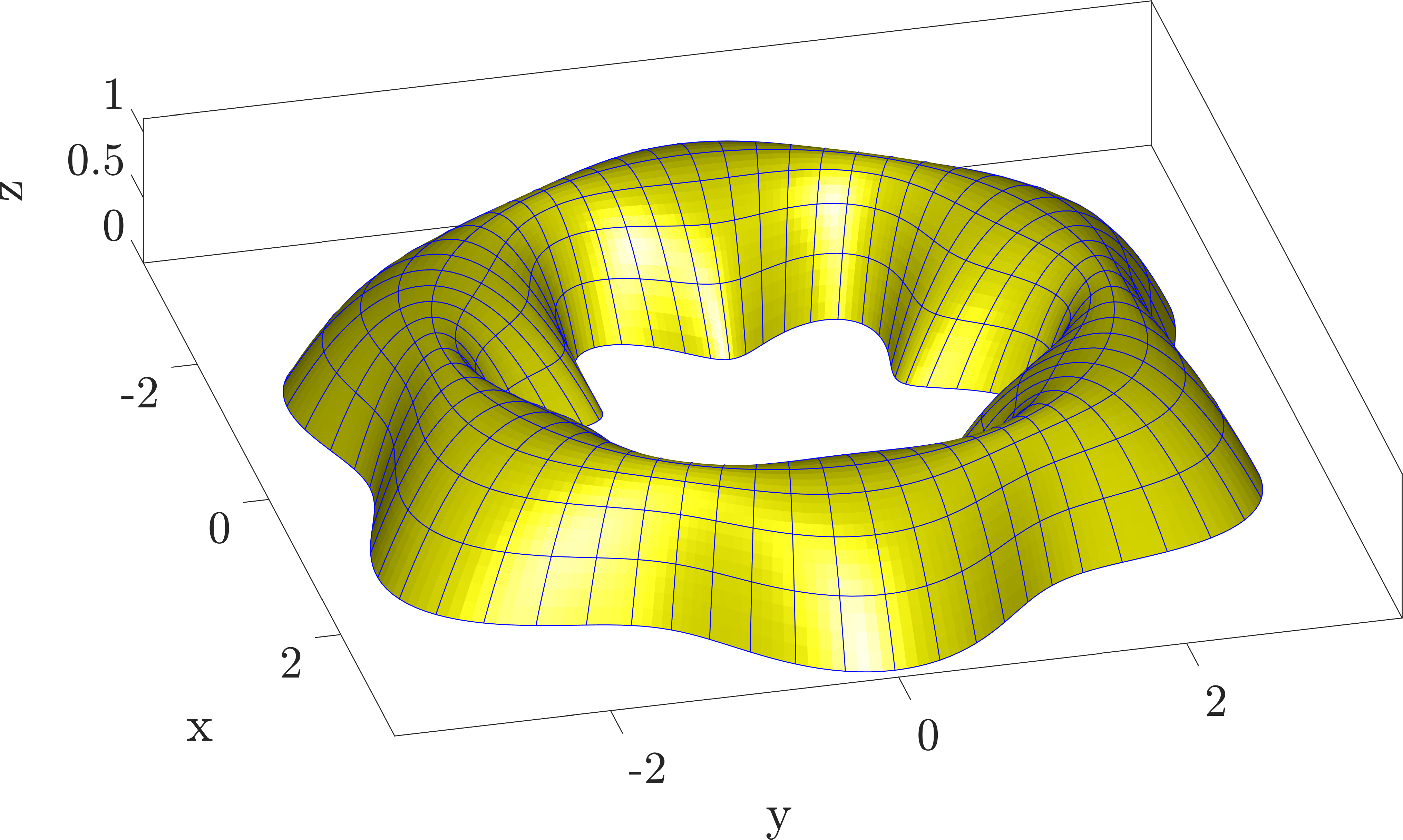}
	\end{minipage}
	\begin{minipage}{.48\textwidth}\scriptsize\flushleft\vspace*{2mm}
		\begin{tabular}[h]{ll}
			Geometry: & Flower shaped shell\\
			& see \autoref{eq:ex3geom} \\
			& $t = \SI{0.1}{}$\\[.25 cm]
			Material parameters: & $E = \SI{10}{}$ \\
			& $\nu = \SI{0.3}{}$\\
			& $\alpha_{\t{s}} = 1.0$ \\[.25 cm]
			Load: & $\vek{f} = [-1\cdot t^3,\,-2\cdot t^3,\,-3\cdot t^3]^\T$\\[.1 cm]
			& $\vek{c} = \vek{0}$\\[.25 cm]
			Support: & Clamped edges at inner and \\
			&outer boundary
		\end{tabular}
	\end{minipage}
	\caption{Definition of flower shaped shell problem.}
	\label{fig:overflower}
\end{figure}

The geometry for the last example is taken from \cite{Schoellhammer_2018a}. The surface is rather complex, but suitable to feature smooth solutions of all physical fields and, consequently, higher-order convergence rates can be achieved. The geometry of the middle surface is given by
\begin{align}
	\vek{x}_\Gamma(r,\,s) &= \begin{bmatrix}
		(A - C) \cos(\theta)\\
		(A - C) \sin(\theta)\\
		1-s^2
	\end{bmatrix}
	\intertext{with:}
	\begin{split}
		&r,\,s \in [-1,\,1]\ ,\ A = 2.3\ ,\ B = 0.8 \\
		&\theta(r) = \pi(r+1) \\
		&C(r,\,s) = s [B + 0.3\cos(6\theta)]\ \label{eq:ex3geom}
	\end{split}
\end{align}
and illustrated in \autoref{fig:overflower}. The shell with the thickness $t = 0.1$ is loaded in all three directions and the material parameters are Young's modulus $E = \SI{10}{}$, Poisson's ratio $\nu = \SI{0.3}{}$. The middle surface features varying principal curvatures and the curved boundaries are clamped.

Following the same rationale as in \cite{Schoellhammer_2018a}, the force equilibrium of \autoref{eq:sff} and moment equilibrium of \autoref{eq:sfm} are computed in strong form and may be called residual errors. In particular, the $L_2$-norms of the residual errors are calculated as
\begin{align}\label{eq:resf}
\varepsilon_{\t{rel,residual,F}}^2 &= \dfrac{\int_\Gamma \left[ \divG{\mat{n}^{\t{real}}_\Gamma} + \mat{Q}\cdot\divG{\mat{q}_\Gamma} + \mat{H}\cdot(\mat{q}_\Gamma\cdot\nG) + \vek{f} \right]^2\ \d A}{\int_\Gamma \vek{f}^2\ \d A} \ ,\\[.25 cm]
\varepsilon_{\t{residual,M}}^2 &= \int_\Gamma \left[\mat{P}\cdot\divG{\mat{m}_ {\Gamma}} - \mat{q}_\Gamma \cdot\nG  + \vek{c}\right]^2\ \d A\ . \label{eq:resm}
\end{align}
The computation of the residual errors requires second-order surface derivatives which is another argument to employ IGA in this work. The theoretical optimal order of convergence in the residual errors is $\mathcal{O}(p-1)$ due to the presence of second order derivatives.\par

\begin{figure}[hbt]
	\centering
	\subfloat[displacement $\vek{u}$]{\includegraphics[width=.36\textwidth, height = .35\textwidth, keepaspectratio]{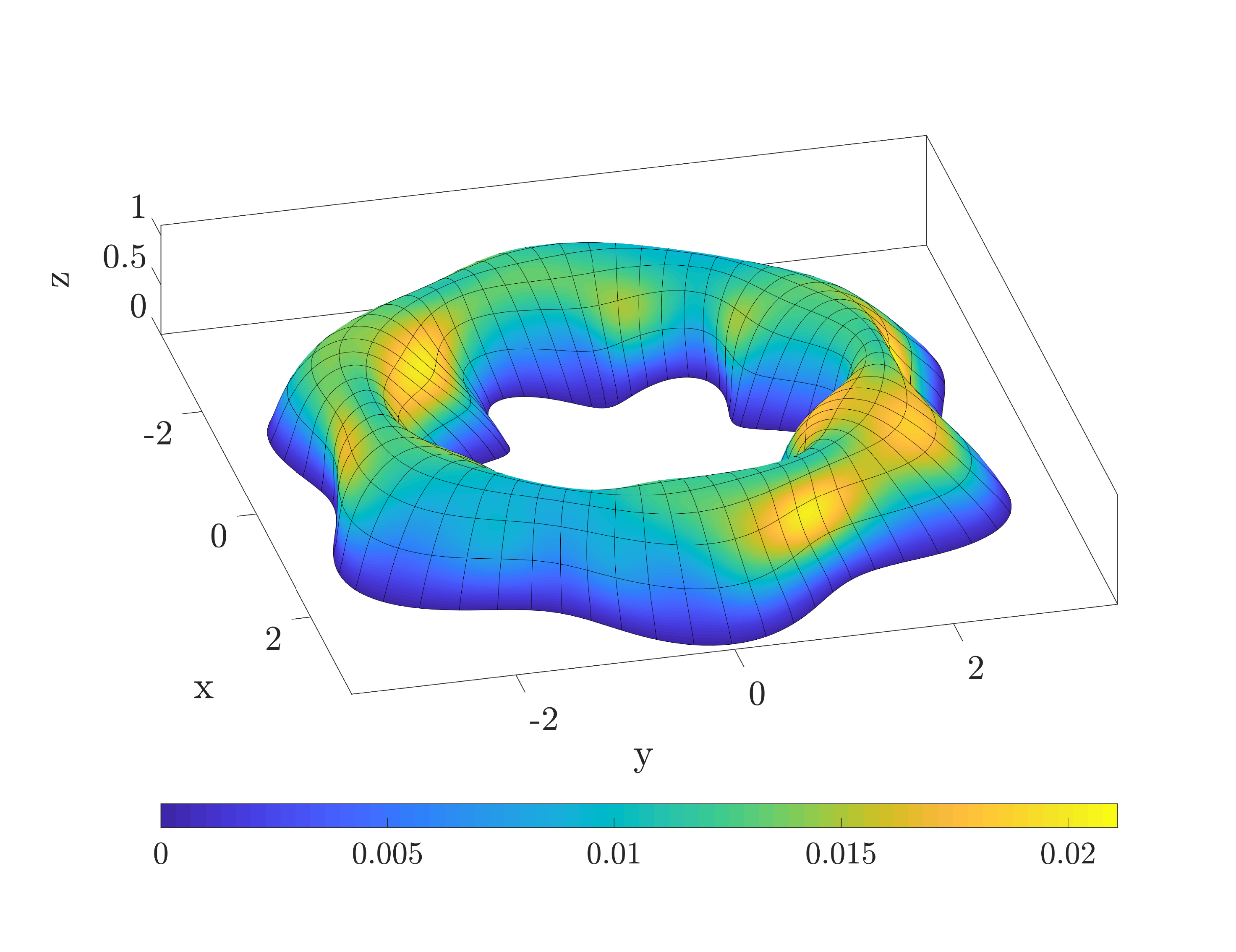}\label{fig:flowershelldisp}}
	\hfil
	\subfloat[force equilibrium]{\includegraphics[width=.3\textwidth, height = .3\textwidth, keepaspectratio]{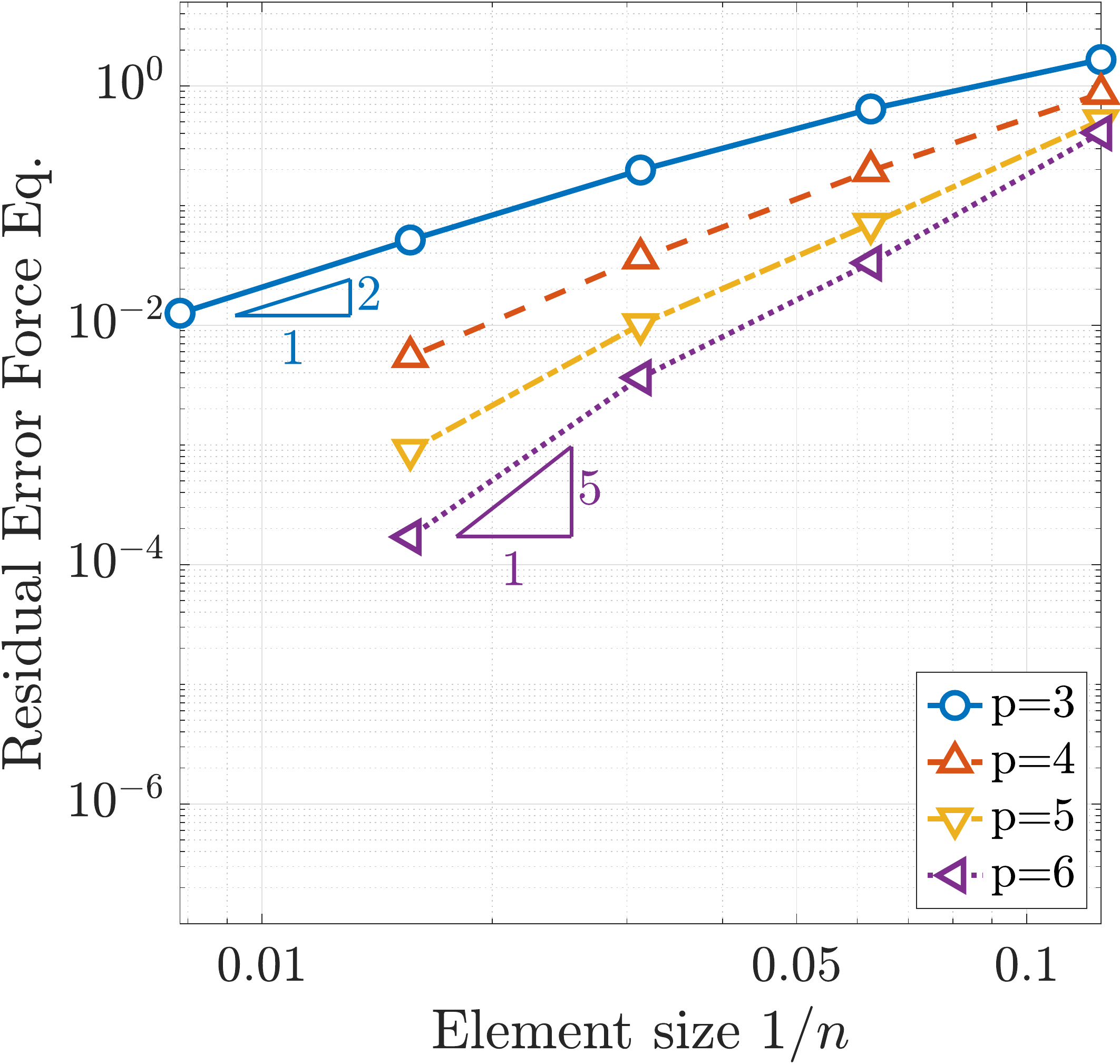}\label{fig:flowershelll2F}}
	\hfil
	\subfloat[moment equilibrium]{\includegraphics[width=.3\textwidth, height = .3\textwidth, keepaspectratio]{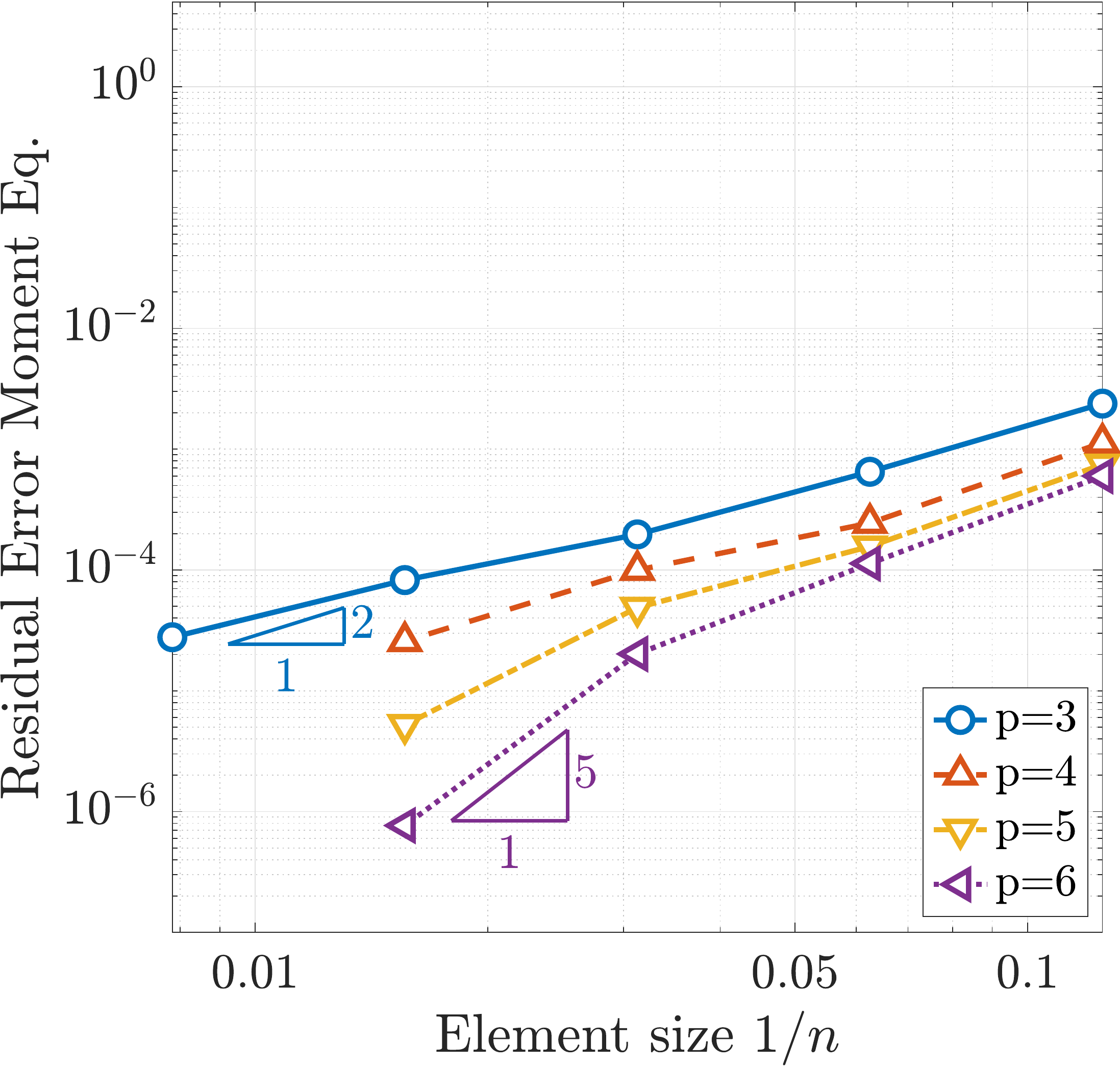}\label{fig:flowershelll2M}}
	\caption{Flower shaped shell: (a) Displacement $\vek{u}$ of flower shaped shell (scaled by one order of magnitude), (b) residual error of the force equilibrium $\varepsilon_{\t{rel,residual,F}}$, (c) residual error of the moment equilibrium $\varepsilon_{\t{residual,M}}$.}
\end{figure}

The numerical solution is presented in \autoref{fig:flowershelldisp} and scaled by one order of magnitude. The results of the convergence study are shown in \autoref{fig:flowershelll2F} and \autoref{fig:flowershelll2M}. The polynomial orders are varied up to $p=6$. Due to the complex geometry and boundary layer effects, the pre-asymptotic range is rather pronounced. Therefore, the results of the coarser levels $n = \lbrace 2,\,4\rbrace$ are omitted in the results. With a sufficiently small knot span size (element size), the expected higher-order convergence rates are achieved in both residual errors \autoref{eq:resf} and \autoref{eq:resm}.

The stored elastic energy at the finest level with a polynomial order $p = 6$ is $\mathfrak{e} = \SI{5.05297916e-04}{}$, which may be seen as an overkill solution. This stored elastic energy may be used for future benchmarking when the computation of the residual errors is not desired.
\section{Conclusions and Outlook}
\label{sec:conc}

The linear Reissner-Mindlin shell theory is reformulated in terms of the TDC using a global Cartesian coordinate system. The total deformation is split into the deformation of the middle surface and the rotation of the normal vector. The latter is expressed with a difference vector. The resulting shell equations are applicable to explicitly (parametrized) or implicitly defined shell geometries. Therefore, the proposed reformulation may seen as a generalization of the classical shell theory which is based on local coordinates resulting from a parametrization of the middle surface.\par

The TDC-based strong form of the force and moment equilibriums is expressed in terms of stress resultants such as membrane forces, bending moments and transverse shear forces using tensor notation. The weak form is consistently derived from the strong form including boundary terms. In the numerical results, the error in the strong form of the shell BVP is integrated over the domain. The resulting residual errors converge with the expected higher orders provided that the involved physical fields are sufficiently smooth. An advantage of investigating the residual errors is that they may be computed for general shell geometries, loadings and boundary conditions without the need for analytical solutions.\par

The weak form is discretized with an isogeometric approach which may be seen as a realization of the Surface FEM using NURBS as trial and test functions. Despite the fact that standard $C^0$-continuous FEM shape functions are sufficient, the high continuity of the isogeometric approach is preferred herein due to the continuous normal vector of the discrete surface, improved convergence properties, and the ability to compute the residual errors. It is clear that in this particular case, the surface mesh, i.e., a NURBS patch, implies a parametrization. This enables a comparison with the classical shell theory and it is seen that the two approaches are equivalent in case of parametrized surfaces, nevertheless, the implementations 
vary significantly. For a general purpose finite element framework, the TDC-based formulation is beneficial, because significant parts of the implementation needed for shells can be shifted to the underlying finite element technology and may be recycled when other PDEs on surfaces are considered.\par

In the numerical results, classical and new benchmark tests are presented and optimal higher-order convergence rates are achieved. There is a large potential in the parametrization-free reformulation of shell models, because the obtained PDEs may be discretized with new finite element techniques such as TraceFEM or CutFEM based on implicitly defined surfaces. In this case, neither the problem statement nor the discretization is based on a parametrization.

\newpage

\clearpage
\newpage

\bibliographystyle{schanz}
\addcontentsline{toc}{section}{\refname}\bibliography{\pathToBibFile}

\begin{thebibliography}{10}

\bibitem{Babuska_1971a}
Babu{\v{s}}ka, I.: Error-bounds for finite element method.
\newblock \emph{Numer. Math.}, \textbf{16}, 322--333, 1971.

\bibitem{Basar_1985a}
Ba{\c{s}}ar, Y.; Kr\"{a}tzig, W.B.: \emph{Mechanik der {F}l\"{a}chentragwerke}.
\newblock Vieweg$+$Teubner Verlag, Braunschweig, 1985.

\bibitem{Bathe_2000a}
Bathe, K.J.; Iosilevich, A.; Chapelle, D.: An evaluation of the {MITC} shell
  elements.
\newblock \emph{Computers \& Structures}, \textbf{75}, 1--30, 2000.

\bibitem{Belytschko_1985a}
Belytschko, T.; Stolarski, H.; Liu, W.K.; Carpenter, N.; Ong, J.S.J: Stress
  projection for membrane and shear locking in shell finite elements.
\newblock \emph{Comp. Methods Appl. Mech. Engrg.}, \textbf{51}, 221--258, 1985.

\bibitem{Benson_2010a}
Benson, D.J.; Bazilevs, Y.; Hsu, M.C.; Hughes, T.J.R.: {I}sogeometric shell
  analysis: {T}he {R}eissner-{M}indlin shell.
\newblock \emph{Comp. Methods Appl. Mech. Engrg.}, \textbf{199}, 276--289,
  2010.

\bibitem{Bieber_2018a}
Bieber, S.; Oesterle, B.; Ramm, E.; Bischoff, M.: A variational method to avoid
  locking--independent of the discretization scheme.
\newblock \emph{Internat. J. Numer. Methods Engrg.}, \textbf{114}, 801--827,
  2018.

\bibitem{Bischoff_2017a}
Bischoff, M.; Ramm, E.; Irslinger., J.: {M}odels and {F}inite {E}lements for
  {T}hin-Walled {S}tructures.
\newblock  \emph{Encyclopedia of Computational Mechanics Second Edition (eds E.
  Stein, R. Borst and T. J. Hughes)}, 2017.

\bibitem{Blaauwendraad_2014a}
Blaauwendraad, J.; Hoefakker, J.H.: \emph{Structural {S}hell {A}nalysis}, Vol.
  200, \emph{Solid Mechanics and Its Applications}.
\newblock {Sprin\-ger}, Berlin, 2014.

\bibitem{Brezzi_1974}
Brezzi, F.: On the existence, uniqueness and approximation of saddle-point
  problems arising from {L}agrange multipliers.
\newblock \emph{RAIRO Anal. Num{\'{e}}r.}, \textbf{R-2}, 129--151, 1974.

\bibitem{Burman_2015a}
Burman, E.; Claus, S.; Hansbo, P.; Larson, M.G.; Massing, A.: Cut{FEM}:
  {D}iscretizing geometry and partial differential equations.
\newblock \emph{Internat. J. Numer. Methods Engrg.}, \textbf{104}, 472--501,
  2015.

\bibitem{Burman_2018a}
Burman, E.; Elfverson, D.; Hansbo, P.; Larson, M.G.; Larsson, K.: Shape
  optimization using the cut finite element method.
\newblock \emph{Comp. Methods Appl. Mech. Engrg.}, \textbf{328}, 242--261,
  2018.

\bibitem{Calladine_1983a}
Calladine, C.~R.: \emph{Theory of {S}hell {S}tructures}.
\newblock Cambridge University Press, Cambridge, 1983.

\bibitem{Cenanovic_2016a}
Cenanovic, M.; Hansbo, P.; Larson, M.G.: Cut finite element modeling of linear
  membranes.
\newblock \emph{Comp. Methods Appl. Mech. Engrg.}, \textbf{310}, 98--111, 2016.

\bibitem{Chapelle_1998a}
Chapelle, D.; Bathe, K.J.: Fundamental considerations for the finite element
  analysis of shell structures.
\newblock \emph{Computers \& Structures}, \textbf{66}, 19--36, 1998.

\bibitem{Chapelle_2000a}
Chapelle, D.; Bathe, K.J.: The mathematical shell model underlying general
  shell elements.
\newblock \emph{Internat. J. Numer. Methods Engrg.}, \textbf{48}, 289--313,
  2000.

\bibitem{Cirak_2001a}
Cirak, F.; Ortiz, M.: Fully {C}1-conforming subdivision elements for finite
  deformation thin-shell analysis.
\newblock \emph{Internat. J. Numer. Methods Engrg.}, \textbf{51}, 813--833,
  2001.

\bibitem{Cirak_2000a}
Cirak, F.; Ortiz, M.; Schr\"oder, P.: Subdivision surfaces: {A} new paradigm
  for thin-shell finite-element analysis.
\newblock \emph{Internat. J. Numer. Methods Engrg.}, \textbf{47}, 2039--2072,
  2000.

\bibitem{Cottrell_2009a}
Cottrell, J.A.; Hughes, T.J.R.; Bazilevs, Y.: \emph{Isogeometric {A}nalysis:
  {T}oward {I}ntegration of {CAD} and {FEA}}.
\newblock {John Wiley \& Sons}, Chichester, 2009.

\bibitem{Delfour_1994a}
Delfour, M.C.; Zol\'{e}sio, J.P.: Shape {A}nalysis via {O}riented {D}istance
  {F}unctions.
\newblock \emph{J. Funct. Anal.}, \textbf{123}, 129--201, 1994.

\bibitem{Delfour_1995a}
Delfour, M.C.; Zol\'{e}sio, J.P.: A {B}oundary {D}ifferential {E}quation for
  {T}hin {S}hells.
\newblock \emph{J. Differential Equations}, \textbf{119}, 426--449, 1995.

\bibitem{Delfour_1996a}
Delfour, M.C.; Zol\'{e}sio, J.P.: Tangential {D}ifferential {E}quations for
  {D}ynamical {T}hin {S}hallow {S}hells.
\newblock \emph{J. Differential Equations}, \textbf{128}, 125--167, 1996.

\bibitem{Delfour_2011a}
Delfour, M.C.; Zol{\'e}sio, J.P.: \emph{{S}hapes and {G}eometries: {M}etrics,
  {A}nalysis, {D}ifferential {C}alculus, and {O}ptimization}.
\newblock SIAM, Philadelphia, 2011.

\bibitem{Demlow_2009a}
Demlow, A.: Higher-order finite element methods and pointwise error estimates
  for elliptic problems on surfaces.
\newblock \emph{SIAM J. Numer. Anal.}, \textbf{47}, 805--827, 2009.

\bibitem{Dornisch_2014a}
Dornisch, W.; Klinkel, S.: Treatment of {R}eissner-{M}indlin shells with kinks
  without the need for drilling rotation stabilization in an isogeometric
  framework.
\newblock \emph{Comp. Methods Appl. Mech. Engrg.}, \textbf{276}, 35--66, 2014.

\bibitem{Dornisch_2013a}
Dornisch, W.; Klinkel, S.; Simeon, B.: Isogeometric {R}eissner-{M}indlin shell
  analysis with exactly calculated director vectors.
\newblock \emph{Comp. Methods Appl. Mech. Engrg.}, \textbf{253}, 491--504,
  2013.

\bibitem{Dziuk_1988a}
Dziuk, G.: \emph{Finite Elements for the {B}eltrami operator on arbitrary
  surfaces}, Chapter~6,  142--155.
\newblock {Sprin\-ger}, Berlin, 1988.

\bibitem{Dziuk_2013a}
Dziuk, G.; Elliott, C.M.: Finite element methods for surface {PDE}s.
\newblock \emph{Acta Numerica}, \textbf{22}, 289--396, 2013.

\bibitem{Echter_2013a}
Echter, R.; Oesterle, B.; Bischoff, M.: A hierarchic family of isogeometric
  shell finite elements.
\newblock \emph{Comp. Methods Appl. Mech. Engrg.}, \textbf{254}, 170--180,
  2013.

\bibitem{Elfverson_2018a}
Elfverson, D.; Larson, M.G.; Larsson, K.: A {N}ew {L}east {S}quares
  {S}tabilized {N}itsche {M}ethod for {C}ut {I}sogeometric {A}nalysis.
\newblock \emph{ArXiv e-prints}, 2018.
\newblock ArXiv: 1804.05654.

\bibitem{Farshad_1992a}
Farshad, M.: \emph{Design and {A}nalysis of {S}hell {S}tructures}.
\newblock {Sprin\-ger}, 1992.

\bibitem{Franca_1988a}
Franca, L.P.; Hughes, T.J.R.: Two classes of mixed finite element methods.
\newblock \emph{Comp. Methods Appl. Mech. Engrg.}, \textbf{69}, 89--129, 1988.

\bibitem{Fries_2018b}
Fries, T.P.: Higher-order surface {FEM} for incompressible {N}avier-Stokes
  flows on manifolds.
\newblock \emph{Int. J. Numer. Methods Fluids}, \textbf{88}, 55--78, 2018.

\bibitem{Fries_2017a}
Fries, T.P.; Omerovi\'{c}, S.; Sch\"{o}llhammer, D.; Steidl, J.: Higher-order
  meshing of implicit geometries - {P}art I: {I}ntegration and interpolation in
  cut elements.
\newblock \emph{Comp. Methods Appl. Mech. Engrg.}, \textbf{313}, 759--784,
  2017.

\bibitem{Fries_2017b}
Fries, T.P.; Sch{\"o}llhammer, D.: Higher-order meshing of implicit geometries
  - {P}art II: {A}pproximations on manifolds.
\newblock \emph{Comp. Methods Appl. Mech. Engrg.}, \textbf{326}, 270--297,
  2017.

\bibitem{Grande_2016a}
Grande, J.; Reusken, A.: A higher order finite element method for partial
  differential equations on surfaces.
\newblock \emph{SIAM}, \textbf{54}, 388--414, 2016.

\bibitem{Gurtin_1975a}
Gurtin, M.E.; Murdoch, I.A.: A continuum theory of elastic material surfaces.
\newblock \emph{Archive for Rational Mechanics and Analysis}, \textbf{57},
  1975.

\bibitem{Hansbo_2014a}
Hansbo, P.; Larson, M.G.: Finite element modeling of a linear membrane shell
  problem using tangential differential calculus.
\newblock \emph{Comp. Methods Appl. Mech. Engrg.}, \textbf{270}, 1--14, 2014.

\bibitem{Hansbo_2017a}
Hansbo, P.; Larson, M.G.: Continuous/discontinuous finite element modelling of
  {K}irchhoff plate structures in {$\mathbb{R}^3$} using tangential
  differential calculus.
\newblock \emph{Comput. Mech.}, \textbf{60}, 693--702, 2017.

\bibitem{Hansbo_2015a}
Hansbo, P.; Larson, M.G.; Larsson, F.: Tangential differential calculus and the
  finite element modeling of a large deformation elastic membrane problem.
\newblock \emph{Comput. Mech.}, \textbf{56}, 87--95, 2015.

\bibitem{Hansbo_2014b}
Hansbo, P.; Larson, M.G.; Larsson, K.: Variational formulation of curved beams
  in global coordinates.
\newblock \emph{Comput. Mech.}, \textbf{53}, 611--623, 2014.

\bibitem{Hughes_2005a}
Hughes, T.J.R.; Cottrell, J.A.; Bazilevs, Y.: Isogeometric analysis: {CAD},
  finite elements, {NURBS}, exact geometry and mesh refinement.
\newblock \emph{Comp. Methods Appl. Mech. Engrg.}, \textbf{194}, 4135--4195,
  2005.

\bibitem{Jankuhn_2017a}
Jankuhn, T.; Olshanskii, M.A.; Reusken, A.: Incompressible {F}luid {P}roblems
  on {E}mbedded {S}urfaces {M}odeling and {V}ariational and {F}ormulations.
\newblock \emph{ArXiv e-prints}, 2017.
\newblock ArXiv: 1702.02989.

\bibitem{Kiendl_2009a}
Kiendl, J.; Bletzinger, K.-U.; Linhard, J.; W{\~{A}}$\sfrac{1}{4}$chner, R.:
  Isogeometric shell analysis with {K}irchhoff-{L}ove elements.
\newblock \emph{Comp. Methods Appl. Mech. Engrg.}, \textbf{198}, 3902--3914,
  2009.

\bibitem{Kiendl_2017a}
Kiendl, J.; Marino, E.; De{\^{A}}{~}Lorenzis, L.: Isogeometric collocation for
  the {R}eissner-{M}indlin shell problem.
\newblock \emph{Comp. Methods Appl. Mech. Engrg.}, \textbf{325}, 645--665,
  2017.

\bibitem{Kirchhoff_1850a}
Kirchhoff, G.: \"{U}ber das {G}leichgewicht und die {B}ewegung einer
  elastischen {S}cheibe.
\newblock \emph{Journal f\"{u}r die reine und angewandte Mathematik (Crelles
  Journal)}, \textbf{40}, 51--88, 1850.

\bibitem{Ko_2017a}
Ko, Y.; Lee, P.S.; Bathe, K.J.: A new {MITC}4+ shell element.
\newblock \emph{Computers \& Structures}, \textbf{182}, 404--418, 2017.

\bibitem{Love_1888a}
Love, A.E.H.: On the small vibrations and deformations of thin elastic shells.
\newblock \emph{Philos. Trans. R. Soc.}, \textbf{179}, 491ff, 1888.

\bibitem{Nguyen-Thanh_2017a}
Nguyen-Thanh, N.; Zhou, K.; Zhuang, X.; Areias, P.; Nguyen-Xuan, H.; Bazilevs,
  Y.; Rabczuk, T.: Isogeometric analysis of large-deformation thin shells using
  {RHT}-splines for multiple-patch coupling.
\newblock \emph{Comp. Methods Appl. Mech. Engrg.}, \textbf{316}, 1157--1178,
  2017.

\bibitem{Oesterle_2016a}
Oesterle, B.; Ramm, E.; Bischoff, M.: A shear deformable, rotation-free
  isogeometric shell formulation.
\newblock \emph{Comp. Methods Appl. Mech. Engrg.}, \textbf{307}, 235--255,
  2016.

\bibitem{Olshanskii_2017a}
Olshanskii, M.A.; Reusken, A.: Trace finite element methods for {PDE}s on
  surfaces.
\newblock \emph{Lecture Notes in Computational Science and Engineering},
  \textbf{121}, 211--258, 2017.

\bibitem{Olshanskii_2017b}
Olshanskii, M.A.; Xu, X.: A trace finite element method for {PDE}s on evolving
  surfaces.
\newblock \emph{SIAM}, \textbf{39}, A1301--A1319, 2017.

\bibitem{Osher_2003a}
Osher, S.; Fedkiw, R.P.: \emph{Level {S}et {M}ethods and {D}ynamic {I}mplicit
  {S}urfaces}.
\newblock {Sprin\-ger}, Berlin, 2003.

\bibitem{Onate_2013a}
O$\tilde{n}$ate, E.: \emph{Structural {A}nalysis with the {F}inite {E}lement
  {M}ethod {L}inear {S}tatics}, Vol. 2. Beams, Plates and Shells, \emph{Lecture
  Notes on Numerical Methods in Engineering and Sciences}.
\newblock {Sprin\-ger}, Berlin, 2013.

\bibitem{Reissner_1945a}
Reissner, E.: The effect of transverse shear deformation on the bending of
  elastic plates.
\newblock \emph{ASME J. Appl. Mech.}, \textbf{12}, A69--77, 1945.

\bibitem{Reusken_2014a}
Reusken, A.: Analysis of trace finite element methods for surface partial
  differential equations.
\newblock \emph{IMA J. Numer. Anal.}, \textbf{35}, 1568--1590, 2014.

\bibitem{Schoellhammer_2018a}
Sch{\"o}llhammer, D.; Fries, T.P.: Kirchhoff-{L}ove shell theory based on
  tangential differential calculus.
\newblock \emph{Comput. Mech.}, 2018.
\newblock DOI: 10.1007/s00466-018-1659-5.

\bibitem{Sethian_1999b}
Sethian, J.A.: \emph{Level {S}et {M}ethods and {F}ast {M}arching Methods}.
\newblock Cambridge University Press, Cambridge, 2nd edition, 1999.

\bibitem{Simo_1989a}
Simo, J.C.; Fox, D.D.: On a stress resultant geometrically exact shell model.
  {P}art I: {F}ormulation and optimal parametrization.
\newblock \emph{Comp. Methods Appl. Mech. Engrg.}, \textbf{72}, 267--304, 1989.

\bibitem{Simo_1989b}
Simo, J.C.; Fox, D.D.; Rifai, M.S.: On a stress resultant geometrically exact
  shell model. {P}art II: {T}he linear theory; {C}omputational aspects.
\newblock \emph{Comp. Methods Appl. Mech. Engrg.}, \textbf{73}, 53--92, 1989.

\bibitem{Tepole_2015a}
Tepole, A.B.; Kabaria, H.; Bletzinger, K.U.; Kuhl, E.: Isogeometric
  {K}irchhoff-{L}ove shell formulations for biological membranes.
\newblock \emph{Comp. Methods Appl. Mech. Engrg.}, \textbf{293}, 328--347,
  2015.

\bibitem{Opstal_2015a}
van Opstal, T.M.; van Brummelen, E.H.; van Zwieten, G.J.: A
  finite-element/boundary-element method for three-dimensional,
  large-displacement fluid-structure-interaction.
\newblock \emph{Comp. Methods Appl. Mech. Engrg.}, \textbf{284}, 637--663,
  2015.

\bibitem{Wempner_2002a}
Wempner, G.; Talaslidis, D.: \emph{{M}echanics of {S}olids and {S}hells:
  {T}heories and {A}pproximations}.
\newblock CRC Press LLC, Florida, 2002.

\bibitem{Yang_2000a}
Yang, H.T.Y.; Saigal, S.; Masud, A.; Kapania, R.K.: A survey of recent shell
  finite elements.
\newblock \emph{Internat. J. Numer. Methods Engrg.}, \textbf{47}, 101--127,
  2000.

\bibitem{Zienkiewicz_2013a}
Zienkiewicz, O.; Taylor, R.; Zhu, J.Z.: \emph{The {F}inite {E}lement {M}ethod:
  {I}ts {B}asis and {F}undamentals: {S}eventh {E}dition}.
\newblock Elsevier LTD, Oxford, 7th edition, 2013.

\bibitem{Zingoni_2018a}
Zingoni, A.: \emph{Shell {S}tructures in {C}ivil and {M}echanical
  {E}ngineering: {T}heory and analysis}.
\newblock ICE Publishing, London, 2nd edition, 2018.

\end{thebibliography}
 
\end{document}